\documentclass[12pt]{article}
\usepackage{lineno,hyperref,latexsym,amsmath, amsfonts, amscd, amssymb, verbatim, setspace,fancyhdr}
\begin{document}

\pagestyle{fancy}
\lhead{N. N. LUAN, J.-C. YAO, AND N. D. YEN}
\rhead{\empty}

\newtheorem{Theorem}{Theorem}[section]
\newtheorem{Proposition}[Theorem]{Proposition}
\newtheorem{Remark}[Theorem]{Remark}
\newtheorem{Lemma}[Theorem]{Lemma}
\newtheorem{Corollary}[Theorem]{Corollary}
\newtheorem{Definition}[Theorem]{Definition}
\newtheorem{Example}[Theorem]{Example}
\renewcommand{\theequation}{\thesection.\arabic{equation}}
\normalsize

\setcounter{equation}{0}

\title{\textbf{On Some Generalized Polyhedral Convex Constructions}}

\medskip
\author{Nguyen Ngoc Luan\footnote{Department of
Mathematics and Informatics, Hanoi National University of Education, 136 Xuan Thuy, Hanoi, Vietnam; email: luannn@hnue.edu.vn.},\ \, Jen-Chih Yao\footnote{Center for General Education, China Medical University, Taichung 40402, Taiwan; email: yaojc@mail.cmu.edu.tw.}, \ \, Nguyen Dong Yen\footnote{Institute of Mathematics, Vietnam Academy of
Science and Technology, 18 Hoang Quoc Viet, Hanoi 10307, Vietnam;
email: ndyen@math.ac.vn.}}

\maketitle
\date{}

\medskip
\begin{quote}
\noindent {\bf Abstract.} Generalized polyhedral convex sets, generalized polyhedral convex functions on locally convex Hausdorff topological vector spaces, and the related constructions such as sum of sets, sum of functions, directional derivative, infimal convolution, normal cone, conjugate function, subdifferential, are studied thoroughly in this paper. Among other things, we show how a generalized polyhedral convex set can be characterized via the finiteness of the number of its faces. In addition, it is proved that the infimal convolution of a generalized polyhedral convex function and a polyhedral convex function is a polyhedral convex function. The obtained results can be applied to scalar optimization problems described by generalized polyhedral convex sets and generalized polyhedral convex functions. 



\end{quote}

\section{Introduction}
\setcounter{equation}{0}  

The concepts of \textit{polyhedral convex set} (pcs) -- also called a \textit{convex polyhedron}, and \textit{generalized polyhedral convex set} (gpcs) -- also called a \textit{generalized convex polyhedron}, stand in the crossroad of several mathematical theories. 

\medskip
First, let us briefly review some basic facts about pcs in a finite-dimensional setting. By definition, \textit{a pcs in a finite-dimensional Euclidean space is the intersection of a finite family of closed half-spaces.} (By convention, the intersection of an empty family of closed half-spaces is the whole space. Therefore, emptyset and the whole space are two special polyhedra.) So, a pcs is the solution set of a system of finitely many inhomogenous linear inequalities. This is the analytical definition of a polyhedral convex set.

\medskip
According to Klee~\cite[Theorem~2.12]{Klee_1959} and Rockafellar~\cite[Theorem~19.1]{Rockafellar_1970}, \textit{for every given convex polyhedron one can find a finite number of points and a finite number of directions such that the polyhedron can be represented as the sum of the convex hull of those points and the convex cone generated by those directions. The converse is also true.} This celebrated theorem, which is a very deep geometrical characterization of pcs, is attributed \cite[p.~427]{Rockafellar_1970} primarily to Minkowski~\cite{Minkowski_1910} and Weyl~\cite{Weyl_1935, Weyl_1953}. By using the result, it is easy to derive fundamental solution existence theorems in linear programming. It is worthy to stress that the above cited representation formula for finite-dimensional pcs has many other applications in mathematics. As an example, we refer to the elegant proofs of the necessary and sufficient second-oder conditions for a local solution and for a locally unique solution in quadratic programming, which were given by Contesse~\cite{Contesse_1980} in 1980; see \cite[pp.~50--63]{Lee_Tam_Yen_2005} for details. 
  
\medskip
For pcs, there is another important characterization: \textit{A closed convex set is a pcs if and only if it has finitely many faces}; see \cite[Theorem~2.12]{Klee_1959} and \cite[Theorem~19.1]{Rockafellar_1970} for details. 
  
\medskip
A bounded pcs is called a \textit{polytope}. Leonhard Euler's Theorem stating \textit{a relation between the numbers of faces of different dimensions of a polytope} is a profound classical result. The reader is referred to \cite[pp.~130--142b]{Grunbaum_2003} for a comprehensive exposition of that theorem and some related results.      

\medskip
Now,  let us discuss the existing facts about pcs and gpcs in an infinite-dimensional setting. According to Bonnans and Shapiro \cite[Definition~2.195]{Bonnans_Shapiro_2000}, a subset of a locally convex Hausdorff topological vector space is said to be a generalized polyhedral convex set (gpcs), or a generalized convex polyhedron, if it is the intersection of finitely many closed half-spaces and a closed affine subspace of that topological vector space. When the affine subspace can be chosen as the whole space, the gpcs is called a polyhedral convex set (pcs), or a convex polyhedron. The theories of generalized linear programming and quadratic programming in \cite[Sections~2.5.7 and 3.4.3]{Bonnans_Shapiro_2000} are mainly based on this concept of gpcs. Some applications of gpcs in Banach spaces can be found in the recent papers by Ban, Mordukhovich and Song~\cite{Ban_Mordukhovich_Song_2011}, Ban and Song~\cite{Ban_Song_2016}, Gfrerer~\cite{Gfrerer_2013, Gfrerer_2014}.

\medskip 
Proposition~2.197 from \cite{Bonnans_Shapiro_2000} tells us that a nonempty gpcs in a locally convex Hausdorff topological vector space has nonempty relative interior. Concerning parametric gpcs in a Banach space, a Hoffman-type lemma was given in \cite[Theorem~2.200]{Bonnans_Shapiro_2000}. Based on that lemma, one can obtain a global Lipschitz continuity property of a generalized polyhedral convex multifunction \cite[Theorem~2.207]{Bonnans_Shapiro_2000}, as well as a local upper Lipschitzian property of a polyhedral multifunction \cite[Theorem~2.208]{Bonnans_Shapiro_2000}.

\medskip
In 2009, using a result related to the Banach open mapping theorem (see, e.g., \cite[Theorem~5.20]{Rudin_1991}), Zheng~\cite[Corollary~2.1]{Zheng_2009} has clarified the relationships between convex polyhedra in Banach spaces and the finite-dimensional convex polyhedra.  

\medskip
Adopting an approach very different from that of Zheng, recently Luan and Yen~\cite{Luan_Yen_2015} have obtained a representation formula for convex polyhedra in locally convex Hausdorff topological vector spaces, which is a comprehensive infinite-dimensional analogue of the above mentioned theorem of Minkowski and Weyl.  In the same paper, the formula has been used for proving solution existence theorems in generalized linear programming and generalized linear vector optimization. Moreover, it allows one to prove \cite[Theorem~4.5]{Luan_Yen_2015} that the weakly efficient solution set of a generalized linear vector optimization problem is the union of finitely many generalized polyhedral convex sets. For the corresponding efficient solution set, a similar result is given in \cite{Luan_2016}, where the relative interior of the dual cone of a polyhedral convex cone is described. Moreover, it is shown that both the solution sets are arcwise connected. Thus the fundamental Arrow-Barankin-Blackwell theorem in linear vector optimization (see \cite{ABB_1953} and \cite{Luc_1989}) has been extended to the locally convex Hausdorff topological vector spaces setting.

\medskip
Functions can be identified with their epigraphs, while sets can be identified with their indicator functions. As explained by Rockafellar \cite[p.~xi]{Rockafellar_1970}, \textit{``These identifications make it easy to pass back and forth between a geometric approach and an analytic approach".} In that spirit, it seems reasonable to call a function \textit{generalized polyhedral convex} when its epigraph is a generalized polyhedral convex set. The introduction of this concept poses an interesting problem. Namely, since the entire Section 19 of~\cite{Rockafellar_1970} is devoted to establishing a variety of basic properties of polyhedral convex sets and polyhedral convex functions (pcf) which have numerous applications afterwards, one may ask whether a similar study can be done for generalized polyhedral convex sets and generalized polyhedral convex functions (gpcf), or not.

\medskip
The aim of the present paper is to solve the above problem. Herein, generalized polyhedral convex sets, generalized polyhedral convex functions on locally convex Hausdorff  topological vector spaces, and the related constructions such as sum of sets,  sum of functions, directional derivative, infimal convolution,  normal cone, conjugate function, subdifferential, will be studied thoroughly. Among other things, we show how a generalized polyhedral convex set can be characterized via the finiteness of the number of its faces. In addition, it is proved that the infimal convolution of a generalized polyhedral convex function and a  polyhedral convex function is a polyhedral convex function. The obtained results can be applied to scalar optimization problems described by generalized polyhedral convex sets and generalized polyhedral convex functions. On one hand, our results can be considered as adequate  extensions of the corresponding classical results in \cite[Section~19]{Rockafellar_1970}. On the other hand, they deepen and develop  the results of~\cite{Luan_2016, Luan_Yen_2015} where only generalized polyhedral convex sets have been considered.

\medskip
Note that Maserick \cite{Maserick_1965} introduced the concept of \textit{convex polytope}, which is very different from the notion of generalized polyhedral convex set in \cite[Definition~2.195]{Bonnans_Shapiro_2000}. On one hand, any convex polytope in the sense of Maserick must have nonempty interior, while a generalized polyhedral convex set in the sense of Bonnans and Shapiro may have empty interior (so it \textit{is not a convex polytope in general}). On the other hand, \textit{there exist convex polytopes in the sense of Maserick which cannot be represented as intersections of finitely many closed half-spaces and a closed affine subspace} of that topological vector space. For example, the closed unit ball $\bar B$ of $c_0$ -- the Banach space of the real sequences $x=(x_1,x_2,\dots)$, $x_i\in\mathbb{R}$ for all $i$, $\lim\limits_{i\to\infty} x_i=0$, with the norm $\|x\|=\sup\{|x_i|\mid i=1,2,\dots\}$ -- is a convex polytope in the sense of Maserick  (see Theorem 4.1 on page 632 in \cite{Maserick_1965}). However, since  $\bar B$ has an infinite number of faces, it cannot be  a generalized polyhedral convex set  in the sense of Bonnans and Shapiro (see Theorem 2.4 in this paper). Subsequently, the concept of convex polytope of \cite{Maserick_1965} has been studied by Maserick and other authors (see, e.g., Durier and Papini \cite{Durier_Papini_1993}, Fonf and Vesely \cite{Fonf_Vesely_2004}). However, after consulting many relevant research works which are available to us, we do hope that the results obtained herein are new. 

\medskip
The organization of the present paper is as follows. Section 2 establishes new facts on generalized polyhedral convex sets. In Section 3, we discuss some basic properties of generalized polyhedral convex functions. Section 4 is devoted to several dual constructions including the concepts of conjugate function and subdifferential of a generalized polyhedral convex function.

\section{Generalized Polyhedral Convex Sets}

We begin this section with a definition of generalized polyhedral convex set.

\medskip
From now on, if not otherwise stated, $X$ is a {\it locally convex Hausdorff topological vector space} (lcHtvs). Denote by $X^*$ the dual space of $X$ and by  $\langle x^*, x \rangle$ the value of~$x^* \in X^*$ at $x \in X$. 
\begin{Definition}{\rm (See \cite[p.~133]{Bonnans_Shapiro_2000})\label{Def_gpcs}
A subset $D \subset X$ is said to be a \textit{generalized polyhedral convex set}, or a \textit{generalized convex polyhedron}, if there exist $x^*_i \in X^*$, $\alpha_i \in \mathbb R$, $i=~1,2,\dots,p$, and a closed affine subspace $L \subset X$, such that 
\begin{equation}\label{eq_def_gpcs}
			D=\big\{ x \in X \mid x \in L,\ \langle x^*_i, x \rangle \leq \alpha_i,\  i=1,\dots,p\big\}.
\end{equation} 
If $D$ can be represented in the form \eqref{eq_def_gpcs} with $L=X$, then we say that it is a \textit{polyhedral convex set}, or a \textit{convex polyhedron}. (Hence, the notion of polyhedral convex set is more specific than that of generalized polyhedral convex set.)}
\end{Definition}

Let $D$ be given as in \eqref{eq_def_gpcs}. According to \cite[Remark~2.196]{Bonnans_Shapiro_2000}, there exists a continuous surjective linear mapping $A$ from $X$ to a lcHtvs $Y$ and a vector $y \in Y$ such that $L=\big\{x \in X \mid Ax=y  \big\}$; then 
\begin{equation}\label{eq_def_gpcs_2}
			D=\big\{ x \in X \mid Ax=y,\  \langle x^*_i, x \rangle \leq \alpha_i,\  i=1,\dots,p\big\}.
\end{equation}   
	
From Definition~\ref{Def_gpcs} it follows that every gpcs is a closed set. If $X$ is finite-dimensional, a subset $D\subset X$ is a gpcs if and only if it is a pcs. For, in that case, we can represent a given affine subspace $L\subset X$ as the solution set of a system of finitely many linear inequalities.
 
\begin{Definition}\label{def_rel_int} {\rm (See \cite[p.~20]{Bonnans_Shapiro_2000})  The \textit{relative interior}  ${\rm ri}C$ of a convex subset $C \subset X$ is the interior of $C$ in the induced topology of the closed affine hull $\overline{\rm aff}C$ of $C$. The \textit{interior} of $C$ is denoted by ${\rm int}C$.}
\end{Definition}

If $X$ is finite-dimensional and $C \subset X$ is a nonempty convex subset, ${\rm ri}C$ is nonempty by \cite[Theorem~6.2]{Rockafellar_1970}. 
If $X$ is infinite-dimensional, it may happen that ${\rm ri}C=\emptyset$ for certain nonempty convex subsets $C \subset X$. To justify the claim, it suffices to choose $X=\ell_2$ -- the Hilbert space of all real sequences $x=(x_k)_{k=1}^{\infty}$ such that $\sum\limits_{k=1}^{\infty} x_k^2 < +\infty$ with the scalar product $\langle x , y \rangle = \sum\limits_{k=1}^{\infty} x_k y_k$. Put $$C=\{x\in \ell_2 \mid x_k \geq 0,\; k=1, 2,\dots\},$$ and observe that ${\rm ri}C={\rm int}C=\emptyset$. If now $C \subset X$ is a nonempty gpcs, it follows, by \cite[Proposition~2.197]{Bonnans_Shapiro_2000}, that ${\rm ri}C \neq \emptyset$. \textit{The latter fact shows that generalized polyhedral convex sets have a nice topological structure.}

\begin{Definition}{\rm (See \cite[p.~162]{Rockafellar_1970})} {\rm A convex subset $F$ of a convex set $C\subset X$ is said to be a \textit{face} of $C$ if for every $x^1, x^2$ in $C$ satisfying $(1-\lambda)x^1+\lambda x^2 \in F$ with $\lambda \in (0,1)$ one has $x^1 \in F$ and $x^2 \in F$. If there exists $x^* \in X^*$ such that $$F=\big\{u \in C \mid \langle x^*, u \rangle =\inf\limits_{x \in C}\langle x^*, x \rangle\big\},$$ then $F$ is called an \textit{exposed face} of $C$. (Therefore, $C$ is not only a face, but also an exposed face of it. The emptyset is a face of $C$, but it is not necessarily an exposed face of~$C$. For example, a nonempty compact convex~$C$ does not have the emptyset as an exposed face of it.)   	
}\end{Definition}

It is necessary to stress that if $F$ is an exposed face of a convex set $C$, then $F$ is a face of $C$. To see that the converse may not true in general, it suffices to choose
\begin{equation*}
C=\left\{(x_1, x_2) \in \mathbb{R}^2 \mid -1 \leq x_1 \leq 1, x_2 \geq -\sqrt{1-x_1^2}\right\}
\end{equation*}
and $F=\{(1,0)\}$.  

\medskip
Theorem~19.1 in \cite{Rockafellar_1970}, which is due to Minkowski \cite{Minkowski_1910} and Weyl \cite{Weyl_1935, Weyl_1953} (see also Klee \cite[Theorem~2.12]{Klee_1959}), is a fundamental result about polyhedral convex sets in finite-dimensional topological vector spaces. In the spirit of that theorem, for a nonempty convex subset $D\subset X$, we are interested in the following properties:
\begin{description}
	\item{\rm (a)} \textit{$D$ is a generalized polyhedral convex set};
	\item{\rm (b)} \textit{There exist $u_1, \dots, u_k \in X$, $v_1, \dots, v_{\ell} \in X$, and a closed linear subspace  $X_0 \subset X$ such that}
	\begin{equation}\label{rep_D}
\begin{aligned}
D=\Bigg\{ \sum\limits_{i=1}^k \lambda_i u_i + \sum\limits_{j=1}^\ell \mu_j v_j \mid &\; \lambda_i \geq 0, \ \forall i=1,\dots,k,    \\ 
&\; \sum\limits_{i=1}^k \lambda_i=1,\ \, \mu_j \geq 0,\ \forall j=1,\dots,\ell   \Bigg\}+X_0;&& 
\end{aligned}
\end{equation}
	\item{\rm (c)} \textit{$D$ is closed and has only a finite number of faces}.
\end{description} 

\medskip
As shown in \cite[Theorem~2.7]{Luan_Yen_2015}, (a) and (b) are equivalent. Now, let us prove that (a) implies (c). 
\begin{Theorem}\label{gcp_has_finitely_many_faces}
	Every generalized polyhedral convex set has only a finite number of faces and all the nonempty faces are exposed.  
\end{Theorem}
\noindent
{\bf Proof.} Let $D$ be a gpcs given by \eqref{eq_def_gpcs}. Set $I=\{1,\dots,p\}$ and $$I(x)=\left\{i \in I \mid \langle x^*_i, x \rangle = \alpha_i \right\}\quad (x \in D).$$ For any subset $J \subset I$, using the definition of face and formula \eqref{eq_def_gpcs}, it is not difficult to show that $F_{J}:=\big\{x \in D \mid \langle x^*_i, x \rangle = \alpha_i, \, \forall i \in J\big\}$ is a face of $D$. 

\medskip
\noindent
{\sc Claim 1.} \textit{Let $x, x' \in D$ and $F$ be a face of $D$. If $x \in F$ and $I(x) \subset I(x')$, then $x' \in F$.}

Indeed, put $x_t:=x-t(x'-x)$ where $t>0$ and observe that $x_t\in L$, because  $x_t=(1+t)x+(-t)x'$ and $x, x'$ belong to the closed affine subspace $L$. For each $i \in I(x) \subset I(x')$, we have
\begin{equation*}
\langle x^*_i,x_t \rangle = \langle x^*_i, x \rangle -t \langle x^*_i, x'-x \rangle=\alpha_i.
\end{equation*}   
 Since $ \langle x^*_j, x \rangle  < \alpha_j$ for all $j \in I \setminus I(x)$, we can find $t>0$ such that $ \langle x^*_j, x_t \rangle  < \alpha_j$ for every $j \in I \setminus I(x)$. Hence, for the chosen $t$, we have $x_t \in D$.  As $x \in F$ and $x=\frac{1}{1+t}x_t + \frac{t}{1+t}x'$,  we must have $x' \in F$. 

\medskip
\noindent
{\sc Claim 2.} \textit{If $F$ is a nonempty face of $D$, then there exists $J \subset I$ such that $F=F_J$. Hence, the number of faces of $D$ is finite. Moreover, $F$ is an exposed face.}

Indeed, given a nonempty face $F$ of $D$, we define $J=\bigcap\limits_{x \in F}I(x)$. It is clear that $F \subset F_J$. To have the inclusion $F_J \subset F$, we select a point $x_0 \in F$ such that the number of elements of $I(x_0)$ is the minimal one among the numbers of elements of $I(x)$, $x \in F$. Let us show that $I(x_0)=J$. Suppose, on the contrary, that $I(x_0) \neq J$. Then there must exist a point $x_1 \in F$ and an index $i_0 \in I(x_0) \setminus I(x_1)$. By the convexity of $F$,  $\bar{x}:=\frac{1}{2}x_0+\frac{1}{2}x_1$ belongs to $F$. Since $\langle x^*_{i_0}, x_1 \rangle < \alpha_{i_0}$, we have $\langle x^*_{i_0}, \bar{x} \rangle < \alpha_{i_0}$, that is, $i_0 \notin I(\bar{x})$. If $j \notin I(x_0)$, i.e., $\langle x^*_j, x_0 \rangle < \alpha_j$, then $\langle x^*_j, \bar{x} \rangle < \alpha_j$; so $j \notin I(\bar{x})$. Thus, $I(\bar{x}) \subset I(x_0)$ and $I(\bar{x}) \neq I(x_0)$. This contradicts the minimality of $I(x_0)$. For any $x \in F_J$, it is clear that $J \subset I(x)$. Since $x_0 \in F$ and $I(x_0)=J \subset I(x)$, by Claim 1 we can assert that $x \in F$. The inclusion $F_J \subset F$ has been proved. Thus $F=F_J$. 

As $J \subset I$ and $I$ is finite, the above obtained result shows that the number of faces of $D$ is finite. 

If $J=\emptyset$, then $F_J=D$. For $x^*:=0$, one has $D={\rm argmin} \big\{\langle x^*, x \rangle \mid x \in D \big\}$; hence~$D$ is an exposed face of it. It follows that $F_\emptyset$ is an exposed face. Now, suppose that $J \neq \emptyset$.  Let $k$ denote the number of elements of $J$. Setting $x^*_J=\frac{1}{k}\sum\limits_{j \in J} (-x^*_j)$, we have $F_J={\rm argmin} \big\{\langle x^*_J, x \rangle \mid x \in D \big\}$. To prove this equality, it suffices to observe that
\begin{equation*}
\langle x^*_J, x \rangle =-\frac{1}{k}\sum\limits_{j \in J}\alpha_j, \ \, \forall x \in F_J. 
\end{equation*} 
and
\begin{equation*}
\langle x^*_J, x \rangle > -\frac{1}{k}\sum\limits_{j \in J}\alpha_j, \ \, \forall x \in D\setminus F_J. 
\end{equation*} 
(The last strict inequality holds because, for any $x \in D\setminus F_J$, there exists $j_0 \in J$ with $\langle - x^*_{j_0}, x \rangle > - \alpha_{j_0}$, while $\langle - x^*_{j}, x \rangle \geq - \alpha_{j}$ for all $j \in J$.) Hence, $F=F_J$ is an exposed face. $\hfill\Box$

\begin{Remark}{\rm 
	The point $x_0$ constructed in the proof of Theorem~\ref{gcp_has_finitely_many_faces} belongs to ${\rm ri}F$. Conversely, for any $\bar{x} \in {\rm ri}F$, $I(\bar{x})$ has the minimality property of $I(x_0)$. The proof of these claims is omitted.} 
\end{Remark}

\begin{Theorem}
	Let $D \subset X$ be a closed convex set with nonempty relative interior. If~$D$ has finitely many faces, then $D$ is a generalized polyhedral convex set.  
\end{Theorem}
\noindent
{\bf Proof.} By our assumption ${\rm ri}D \neq \emptyset$. We first consider the case in where ${\rm int}D \neq \emptyset$. We have $D={\rm int}D \cup \partial D$, where $\partial D=D \setminus {\rm int}D$ is the boundary of $D$. If $\partial D =\emptyset$ then $D=X$ because $D$ is both open and closed in $X$, which is a connected topological space. So $D$ is a convex polyhedron. If $\partial D \neq \emptyset$, we pick a point $\bar{x} \in \partial D$. As $\{\bar{x}\} \cap {\rm int}D=~\emptyset$ and since $\{\bar{x}\}$ and ${\rm int}D$ are convex sets, by the separation theorem \cite[Theorem~3.4 (a)]{Rudin_1991}, there exists $\varphi_{\bar{x}} \in X^* \setminus\{0\}$ such that $\langle \varphi_{\bar{x}}, \bar{x} \rangle \geq  \langle \varphi_{\bar{x}},x \rangle$ for all $x \in {\rm int}D$. Since $D$ is convex and ${\rm int}D \neq \emptyset$, it follows that
\begin{equation}\label{proof_face_1}
\langle \varphi_{\bar{x}}, \bar{x} \rangle \geq  \langle \varphi_{\bar{x}},x \rangle, \, \ \forall x \in D. 
\end{equation} 
Let $\alpha_{\bar{x}}:=\langle \varphi_{\bar{x}}, \bar{x} \rangle$ and $F_{\bar{x}, \varphi_{\bar{x}}}:=\{ x\in D \mid \langle \varphi_{\bar{x}}, x \rangle =\alpha_{\bar{x}} \}$. 
It is easy to show that $F_{\bar{x}, \varphi_{\bar{x}}}$ is a face of $D$ and $\bar{x} \in F_{\bar{x}, \varphi_{\bar{x}}}$. As $D$ has finitely many faces, we can find a finite sequence of points $x_1, \dots, x_k$ in $\partial D$ such that, for every $u \in  \partial D$, there exists $i \in \{1,\dots,k\}$ with $F_{u, \varphi_{u}}=F_{x_i, \varphi_{x_i}}$. Let 
\begin{equation}\label{proof_face_2}
D':=\big\{x \in X \mid \langle \varphi_{x_i}, x \rangle \leq \alpha_{x_i}, \, i=1\,\dots,k \big\}.
\end{equation} 
By the construction of $\varphi_{x_i}$, $i=1,\dots,k$, and by \eqref{proof_face_1}, we have $D \subset D'$.  To show that $D'=D$, suppose the contrary: There exists $u_1 \in D'\setminus D$. Select a point $u_0 \in {\rm int}D$. Let $[u_0, u_1]:=\big\{(1-t)u_0+tu_1 \mid t \in [0,1]\big\}$ denote the segment joining $u_0$ and $u_1$. Since $[u_0, u_1] \cap D$ is a nonempty closed convex set, $T:=\{t \in [0,1] \mid u_t:=(1-t)u_0+tu_1 \in D\}$ is a closed convex subset of $[0,1]$. Note that $0 \in T$, but $1 \notin T$. Hence, $T=[0, \bar{t}]$ for some $\bar{t} \in [0,1)$. As $u_0 \in {\rm int}D$, we must have $\bar{t} >0$. It is easy to show that $\bar{u}:=\left(1-\bar{t}\right)u_0+\bar{t}u_1$ belongs to $\partial D$. Hence, $F_{\bar{u}, \varphi_{\bar{u}} }=F_{x_i, \varphi_{x_i}}$ for some $i \in \{1,\dots,k\}$. Since $u_0 \in {\rm int}D$ and $\varphi_{x_i} \neq 0$, from \eqref{proof_face_1} it follows that 
\begin{equation}\label{proof_face_3}
\langle \varphi_{x_i}, u_0 \rangle < \alpha_{x_i}.
\end{equation}
As $\bar{u} \in F_{x_i, \varphi_{x_i}}$, one has 
\begin{equation}\label{proof_face_4}
\langle \varphi_{x_i}, \bar{u} \rangle = \alpha_{x_i}.
\end{equation}
From the equality $\bar{u}=\left(1-\bar{t}\right)u_0+\bar{t}u_1$ we can deduce that $u_1=\frac{1}{\bar{t}}\bar{u}+\left(1-\frac{1}{\bar{t}}\right)u_0$. Since $1-\frac{1}{\bar{t}} < 0$, by \eqref{proof_face_3} and \eqref{proof_face_4} we have
\begin{equation*}
\begin{aligned}
\langle \varphi_{x_i}, u_1 \rangle &=\frac{1}{\bar{t}} \langle \varphi_{x_i}, \bar{u}\rangle + \left(1-\frac{1}{\bar{t}}\right) \langle \varphi_{x_i}, u_0 \rangle\\
&>\frac{1}{\bar{t}} \alpha_{x_i}+\left(1-\frac{1}{\bar{t}}\right)\alpha_{x_i}=\alpha_{x_i}.
\end{aligned}
\end{equation*}
Then we obtain $ \langle \varphi_{x_i}, u_1 \rangle  > \alpha_{x_i}$, contradicting the assumption $u_1 \in D'$. We have thus proved that $D'=D$. Therefore, by \eqref{proof_face_2} we can conclude that $D$ is a polyhedral convex set.

Now, let us consider the case ${\rm int}D=\emptyset$. As ${\rm ri}D \neq \emptyset$, the interior of $D$ in the induced topology of $\overline{\rm aff}D$ is nonempty. Take any $x_0 \in D$. Applying the above result for the closed convex subset $D_0:=D-x_0$ of the lcHtvs $X_0:=\overline{\rm aff}D-x_0$, we find $x^*_i \in X_0^*$ and $\alpha_i \in \mathbb{R}$, $i=1,\dots,m$, such that
\begin{equation}\label{proof_face_5}
D_0=\big\{x \in X_0 \mid \langle x^*_i, x \rangle \leq \alpha_i, \ \, i=1,\dots,m\big\}.
\end{equation} 
By the well-known extension theorem \cite[Theorem~3.6]{Rudin_1991}, we can find $\widetilde{x}^*_i \in X^*$, $i=1,\dots,m$, such that $\langle \widetilde{x}^*_i, x \rangle = \langle x^*_i, x \rangle$ for all $x \in X_0$. Then from \eqref{proof_face_5} it follows~that 
\begin{equation*}
D_0=\big\{x \in X_0 \mid \langle \widetilde{x}^*_i, x \rangle \leq \alpha_i, \ \, i=1,\dots,m\big\}.
\end{equation*}  
As $D=D_0+x_0$, this implies that  
\begin{equation*}
D=\big\{u \in X_0+x_0 \mid \langle \widetilde{x}^*_i, u \rangle \leq \alpha_i+ \langle \widetilde{x}^*_i, x_0 \rangle, \ \, i=1,\dots,m\big\}.
\end{equation*} 
Thus $D$ is a generalized polyhedral convex set. $\hfill\Box$

\medskip
Let us consider the following question: \textit{Whether the image of a generalized polyhedral convex set via a linear mapping from $X$ to $Y$, which are locally convex Hausdorff topological vector spaces, is a generalized polyhedral convex set, or not?} The answers in the affirmative are given in \cite[Theorem~19.3]{Rockafellar_1970} for the case where $X$ and $Y$ are finite-dimensional, in \cite[Lemma~3.2]{Zheng_Yang_2008} for the case where $X$ is a Banach space and~$Y$ is finite-dimensional, and in \cite[Proposition~2.1]{Luan_2016} for the case where $X$ is a lcHtvs and $Y$ is finite-dimensional. When $Y$ is infinite-dimensional, the image may not be a gpcs (see \cite[Example~2.1]{Luan_2016}). In the above mentioned example, one sees that the image of a closed linear subspace of $X$ via a continuous surjective linear operator may be non-closed; hence it cannot be a generalized polyhedral convex set. 

\medskip
The above results motivate the following proposition.
\begin{Proposition}\label{image_gpc}
	Suppose that $T: X \rightarrow Y$ is a linear mapping between locally convex Hausdorff topological vector spaces and $D \subset X$, $Q \subset Y$ are nonempty generalized polyhedral convex sets. Then, $\overline{T(D)}$ is a generalized polyhedral convex set. If $T$ is continuous, then $T^{-1}(Q)$ is a generalized polyhedral convex set.	
\end{Proposition}
\noindent
{\bf Proof.} Suppose that $D$ is of the form \eqref{rep_D}. Then  $T(D)=D'+T(X_0)$, where~$D':={\rm conv}\big\{Tu_i \mid i=1,\dots,k\big\} + {\rm cone}\big\{Tv_j \mid i=1,\dots,\ell\big\}$ with ${\rm conv}\Omega$ denoting the convex hull of a subset $\Omega\subset Y$ and ${\rm cone}M$ denoting the convex cone generated by a subset $M\subset Y$. Since $T(X_0) \subset Y$ is a linear subspace, $\overline{T(X_0)}$ is a closed linear subspace of~$Y$ by \cite[Theorem~1.13~(c)]{Rudin_1991}; so $D'+\overline{T(X_0)}$ is a gpcs by \cite[Theorem~2.7]{Luan_Yen_2015}. In particular, $D'+\overline{T(X_0)}$ is closed. Hence, the inclusion $T(D) \subset D'+\overline{T(X_0)}$ yields 
\begin{equation}\label{proof_image_eq_1}
\overline{T(D)} \subset D'+\overline{T(X_0)}.
\end{equation}
According to \cite[Theorem~1.13 (b)]{Rudin_1991}, we have
\begin{equation}\label{proof_image_eq_2}
D'+\overline{T(X_0)} \subset \overline{T(D)}.
\end{equation}
Combining \eqref{proof_image_eq_1} with \eqref{proof_image_eq_2} implies that $\overline{T(D)}=D'+\overline{T(X_0)}$. Therefore $\overline{T(D)}$ is a generalized polyhedral convex set. 

Now, suppose that $Q \subset Y$ is a gpcs given by $$Q=\big\{y \in Y \mid By=z,\  \langle y^*_j, y \rangle \leq \beta_j,\ j=1,\dots,q\big\},$$ where $B: Y \to Z$ is a continuous linear mapping between two lcHtvs, $z \in Z$ and $y^*_j \in Y^*$, $\beta_j \in \mathbb{R}$, $j=1,\dots,q$. Then we have 
\begin{equation*}
\begin{aligned}
T^{-1}(Q)=&\big\{x \in X \mid B(Tx)=z,\ \langle y^*_j, Tx \rangle \leq \beta_j,\ j=1,\dots,q\big\}\\
=&\big\{x \in X \mid (B\circ T)x=z,\ \langle T^*y^*_j, x \rangle \leq \beta_j, \ j=1,\dots,q\big\},
\end{aligned}
\end{equation*}
where $T^*: Y^* \rightarrow X^*$ is the adjoint operator of $T$. Since $T: X \to Y$ and $B: Y \to Z$ are linear continuous mappings, $B\circ T: X \to Z$ is a continuous linear mapping. Hence, the above expression for $T^{-1}(Q)$ shows that the set is generalized polyhedral convex. 
$\hfill\Box$

\begin{Corollary}\label{sum_gcps}
	If $D_1, \dots, D_m$ are nonempty generalized polyhedral convex sets in $X$, so is $\overline{D_1+\dots+D_m}$.  
\end{Corollary}
\noindent
{\bf Proof.} Consider the linear mapping $T: X^m \to X$ given by $$T(x_1, \dots, x_m)=x_1+\dots+x_m\quad \forall (x_1, \dots, x_m) \in X^m,$$ and observe that $T(D_1 \times \dots \times D_m)=D_1+\dots+D_m$. 
Since $D_k$ is a gpcs in $X$ for $k=1,\dots,m$, 
using Definition \ref{Def_gpcs} one can show that $D_1 \times \dots \times D_m$ is a gpcs in $X^m$. Then, $\overline{T(D_1 \times \dots \times D_m)}$ is a gpcs by Proposition \ref{image_gpc}. Hence, $\overline{D_1+\dots+D_m}$ is a gpcs in $X$. 
$\hfill\Box$

\begin{Remark}\label{Rem_sum_closed_spaces}
	{\rm One may ask: \textit{Whether the closure sign can be removed from Corollary~\ref{sum_gcps}, or not?} When $X$ is a finite-dimensional space, the sum of finitely many pcs in~$X$ is a pcs (see, e.g., \cite[Corollary~2.16]{Klee_1959}, \cite[Corollary~19.3.2]{Rockafellar_1970}). However, when $X$ is an infinite-dimensional space, the sum of a finite family of gpcs may not be a gpcs. To see this, one can choose a suitable space $X$ and closed linear subspaces $X_1, X_2$ of~$X$ so that  $\overline{X_1+X_2}=X$ and $X_1+X_2 \neq X$ (see \cite[Example~3.34]{Bauschke_Combettes_2011} for an example of subspaces in any infinite-dimensional Hilbert space, \cite[Exercise~1.14]{Brezis_2011} for an example in~$\ell^1$, and \cite[Exercise~20, p.~40]{Rudin_1991} for an example in $L^2(-\pi, \pi)$). Clearly, $X_1, X_2$ are gpcs in~$X$. Since $X_1+X_2$ is non-closed, it cannot be a gpcs. 
}\end{Remark}

Concerning the question stated in Remark \ref{Rem_sum_closed_spaces}, in the two following propositions we shall describe some situations where the closure sign can be dropped.   	

\begin{Proposition}\label{sum_pc_gpc2}
	If $D_1, D_2$ are generalized polyhedral convex sets of $X$ and ${\rm aff}D_1$ is finite-dimensional,  then $D_1+D_2$ is a generalized polyhedral convex set.
\end{Proposition}
\noindent
{\bf Proof.} According to \cite[Theorem~2.7]{Luan_Yen_2015}, for each $m\in \{1,2\}$, we can represent $D_m$ as $D_m=D_m'+X_{m,0}$ with $X_{m,0}$ being a closed linear subspace of $X$, $$D_m'={\rm conv} \left\{u_{m,1}, \dots,u_{m,k_m}\right\}+{\rm cone} \left\{v_{m,1},\dots,v_{m, \ell_m}\right\}$$ for some
$u_{m,1},\dots,u_{m,k_m}$, $v_{m,1},\dots,v_{m,\ell_m}$ in $X$. Since ${\rm aff}D_1$ is finite-dimensional, we must have ${\rm dim}X_{1,0} < \infty$. By \cite[Theorem~1.42]{Rudin_1991}, $X_{1,0}+X_{2,0}$ is a closed linear subspace of $X$.  Let $W$ be the finite-dimensional linear subspace generated by the vectors $u_{m,1},\dots,u_{m,k_m}$, $v_{m,1},\dots,v_{m,\ell_m}$, for $m=1,2$. Since $D_1'$ and $D_2'$ are polyhedral convex sets in $W$ due to \cite[Theorem~19.1]{Rockafellar_1970}, $D_1'+D_2'$ is a pcs in $W$ by \cite[Corollary~19.3.2]{Rockafellar_1970}. On account of \cite[Theorem~19.1]{Rockafellar_1970}, one can choose $u_1,\dots,u_k$ in $W$, $v_1,\dots,v_{\ell}$ in $W$ such that $D_1'+D_2'={\rm conv} \left\{u_{i} \mid i=1,\dots,k\right\}+{\rm cone} \left\{v_{j} \mid j=1,\dots,\ell\right\}$. It follows that
\begin{equation*}\label{proof_rep_sumD1D2}
D_1+D_2 = {\rm conv} \left\{u_{i} \mid i=1,\dots,k\right\}+{\rm cone} \left\{v_{j} \mid j=1,\dots,\ell\right\}+X_{1,0}+X_{2,0}.
\end{equation*}  
Recalling that the linear subspace $X_{1,0}+X_{2,0}$ is closed,  we can use \cite[Theorem~2.7]{Luan_Yen_2015} to assert that $D_1+D_2$ is a gpcs. $\hfill\Box$

\begin{Proposition}\label{sum_pc_gpc}
	If $D_1 \subset X$ is a polyhedral convex set and $D_2 \subset X$ is a generalized polyhedral convex set,  then $D_1+D_2$ is a polyhedral convex set.
\end{Proposition}

The proof of this result is based on the next two lemmas.
 \begin{Lemma}\label{rep_pcs_lemma}
	 A nonempty subset $D \subset X$ is polyhedral convex if and only if $D$ admits a representation of the form \eqref{rep_D}, where $X_0$ is a closed linear subspace of finite codimension.
\end{Lemma}
\noindent
{\bf Proof.} Similar to the proof of \cite[Theorem~2.7]{Luan_Yen_2015}. $\hfill\Box$

\begin{Lemma}\label{closed_sum_two_subspaces}
	If $X_1$ and $X_2$ are linear subspaces of $X$ with $X_1$ being closed and finite-codimensional, then $X_1+X_2$ is closed and ${\rm codim}(X_1+X_2) < \infty$.    
\end{Lemma}
\noindent
{\bf Proof.} Since $X_1 \subset X$ is finite-codimensional, there exists a finite-dimensional linear subspace $X_1' \subset X$ such that $X=X_1 \cup X_1'$ and $X_1 \cap X_1'=\{0\}$. Let $\pi_1: X \rightarrow X/X_1$, $\pi_1(x)= x+X_1$ for every $x \in X$, be the canonical projection from $X$ on the quotient space $X/X_1$. It is clear that the operator $\Phi_1: X/X_1 \rightarrow X_1'$, $x' +X_1 \mapsto x'$ for all $x' \in X_1'$, is a linear bijective mapping. On one hand, by \cite[Theorem~1.41({\it a})]{Rudin_1991}, $\pi_1$ is a linear continuous mapping. On the other hand, $\Phi_1$ is a homeomorphism by \cite[Lemma~2.5]{Luan_Yen_2015}. So, the operator $\pi:=\Phi_1\circ \pi_1: X \rightarrow X_1'$ is linear and continuous. Note that $\pi(X_2)$ is closed, because it is a linear subspace of $X_1'$, which is finite-dimensional. Since $\pi$ is continuous and $X_1+X_2=\pi^{-1}\left(\pi(X_2)\right)$, we see that $X_1+X_2$ is closed. The~${\rm codim}X_1 < \infty$ clearly forces ${\rm codim}(X_1+X_2) < \infty$.    
$\hfill\Box$

\medskip
\noindent
{\bf Proof of Proposition~\ref{sum_pc_gpc}.} By Lemma~\ref{rep_pcs_lemma}, there exist $u_{1,1},\dots,u_{1,k_1}$ in $X$,  $v_{1,1},\dots,v_{1,\ell_1}$ in $X$ and a closed finite-codimensional linear subspace $X_{1,0} \subset X$ such that $D_1=D_1'+X_{1,0}$ with $D_1'={\rm conv}\left\{u_{1,1}, \dots,u_{1,k_1}\right\}+{\rm cone} \left\{v_{1,1}, \dots, v_{1,\ell_1}\right\}$. According to \cite[Theorem~2.7]{Luan_Yen_2015}, there exist $u_{2,1},\dots,u_{2,k_2}$ in $X$, $v_{2,1},\dots,v_{2,\ell_2}$ in $X$ and a closed linear subspace $X_{2,0} \subset X$ satisfying $D_2=D_2'+X_{2,0}$ with 
$D_2'={\rm conv} \left\{u_{2,1}, \dots,u_{2,k_2}\right\}+{\rm cone} \left\{v_{2,1}, \dots, v_{2,\ell_2}\right\}.$
 Let $W$ be the finite-dimensional linear subspace generated by the vectors $u_{1,1},\dots,u_{1,k_1}$, $v_{1,1},\dots,v_{1,\ell_1}$, $u_{2,1},\dots,u_{2,k_2}$, $v_{2,1},\dots,v_{2,\ell_2}$. Since $D_1'$ and $D_2'$ are pcs in $W$ by \cite[Theorem~19.1]{Rockafellar_1970}, Corollary~19.3.2 of \cite{Rockafellar_1970} implies that $D_1'+D_2'$ is a pcs. Applying \cite[Theorem~19.1]{Rockafellar_1970} for the pcs $D_1'+D_2'$ of $W$, one can find $u_1,\dots,u_k$ and $v_1,\dots,v_{\ell}$ in $W$ such that $D_1'+D_2'={\rm conv} \left\{u_{i} \mid i=1,\dots,k\right\}+{\rm cone} \left\{v_{j} \mid j=1,\dots,\ell\right\}$. Thus, 
\begin{equation}\label{proof_rep_sumD1D2}
D_1+D_2 = {\rm conv} \left\{u_{i} \mid i=1,\dots,k\right\}+{\rm cone} \left\{v_{j} \mid j=1,\dots,\ell\right\}+X_{1,0}+X_{2,0}.
\end{equation}  
In accordance with Lemma \ref{closed_sum_two_subspaces}, $X_{1,0}+X_{2,0}$ is a closed finite-codimensional linear subspace. Hence, by Lemma \ref{rep_pcs_lemma} and formula \eqref{proof_rep_sumD1D2} we conclude that $D_1+D_2$ is a pcs.  $\hfill\Box$

\medskip
The next result is an extension of \cite[Corollary~19.3.2]{Rockafellar_1970} to an infinite-dimensional setting.
\begin{Corollary}\label{strong_sep_pc_pcs}
	Suppose that $D_1 \subset X$ is a polyhedral convex set and $D_2 \subset X$ is a generalized polyhedral convex set. If $D_1 \cap D_2=\emptyset$, then there exists $x^* \in X^*$ such that
	\begin{equation}\label{strongly_sep_pcs}
	\sup\{\langle x^*, u \rangle \mid u \in D_1 \} < \inf\{\langle x^*, v \rangle \mid v \in D_2 \}.    
	\end{equation}
\end{Corollary}
\noindent
{\bf Proof.} By Proposition~\ref{sum_pc_gpc}, $D_2-D_1=D_2+(-D_1)$ is a polyhedral convex set in~$X$; hence it is closed. Since $D_2-D_1$ is a closed convex set and $0 \notin D_2-D_1$, by the strongly separation theorem \cite[Theorem~3.4 (b)]{Rudin_1991} there exist $x^* \in X^*$ and $\gamma \in \mathbb{R}$ such that 
\begin{equation*}
\langle x^*, 0 \rangle  < \gamma   \leq \langle x^*, x \rangle,  \quad  \forall x \in D_2-D_1. 
\end{equation*}    
This implies that $\sup\{\langle x^*, u \rangle \mid u \in D_1 \}+\gamma \leq \inf\{\langle x^*, v \rangle \mid v \in D_2 \}$; hence the strict inequality \eqref{strongly_sep_pcs} is valid. $\hfill\Box$

\medskip
The assertion of Corollary~\ref{strong_sep_pc_pcs} would be false if $D_1$ is only assumed to be a gpcs. Indeed, an answer in the negative for the question in \cite[Exercise~1.14]{Brezis_2011} assures us that there exist closed affine subspaces  $D_1$ and $D_2$ in $X=\ell^1$ such that one cannot find any $x^* \in X^*\setminus\{0\}$ satisfying  $\sup\{\langle x^*, u \rangle \mid u \in D_1 \}\leq \inf\{\langle x^*, v \rangle \mid v \in D_2 \}$. So, with the chosen generalized polyhedral convex sets $D_1$ and $D_2$, one cannot have \eqref{strongly_sep_pcs} for any $x^*\in X^*=\ell^\infty$. 

\medskip
As in \cite[p.~61]{Rockafellar_1970}, the {\it recession cone} $0^{+}C$ of a convex set $C \subset X$ is given by 
\begin{equation*}
0^{+}C=\big\{v \in X \mid x+tv \in C, \ \forall x \in X,\ \forall t \geq 0\big\}.
\end{equation*} 
If $C$ is nonempty and closed, then $0^+C$ is a closed convex cone, and $v \in X$ belongs to~$0^{+}C$ if and only if there exists $x \in C$ such that $x + tv \in C$ for all $t \geq 0$. These facts are well known \cite[Theorems~8.2 and 8.3]{Rockafellar_1970} for closed convex sets in $\mathbb R^n$. For the general case where $X$ is a lcHtvs, the  facts can be found in \cite[p.~33]{Bonnans_Shapiro_2000}. 

\medskip
We are now in a position to extend Theorem~19.6 from the book of Rockafellar~\cite{Rockafellar_1970}, which was given in $\mathbb{R}^n$, to the case of generalized polyhedral convex in lcHtvs.
\begin{Theorem}\label{cup_gcps}
Suppose that $D_1, \dots, D_m$ are generalized polyhedral convex sets in $X$. Let $D$ be the smallest closed convex subset of $X$ that contains $D_i$ for all $i=1,\dots,m$. Then $D$ is a generalized polyhedral convex set. If at least one of the sets $D_1, \dots, D_m$ is polyhedral convex, then $D$ is a  polyhedral convex set.
\end{Theorem}
\noindent
{\bf Proof.} By removing all the empty sets from the system $D_1, \dots, D_m$, we may assume that $D_i\neq\emptyset$ for all $i\in I:=\{1,\dots,m\}$. Due to \cite[Theorem~2.7]{Luan_Yen_2015}, for each $i\in I$, one can find $u_{i,1},\dots,u_{i,k_i}$ and $v_{i,1},\dots,v_{i,\ell_i}$ in $X$ and a closed linear subspace $X_{i,0} \subset X$ such that
\begin{equation}\label{eq_rep_D_i}
D_i={\rm conv}\{u_{i,1},\dots,u_{i,k_i}\}+{\rm cone}\{v_{i,1},\dots,v_{i,\ell_i}\}+X_{i,0}.
\end{equation}
Since $X_{1,0}+\dots+X_{m,0} \subset X$ is a linear subspace, $X_0:=\overline{X_{1,0}+\dots+X_{m,0}}$ is a closed linear subspace of $X$ by \cite[Theorem~1.13 (c)]{Rudin_1991}. Let
\begin{equation}\label{eq_cup_gcps}
D':={\rm conv}\{u_{i,j}  \mid  i \in I,  j=1,\dots,k_i \} + {\rm cone}\{v_{i,j}  \mid i \in I,  j=1,\dots,\ell_i \} + X_{0}. 
\end{equation} 
On account of \cite[Theorem~2.7]{Luan_Yen_2015}, $D'$ is a gpcs. In particular, $D'$ is convex and closed. From \eqref{eq_rep_D_i} and \eqref{eq_cup_gcps} it follows that $D_i \subset D'$ for every $i \in I$. Hence, by the definition of $D$, we must have $D \subset D'$. Let us show that $D'\subset D$. Since $u_{i, j}$ belongs to $D_i \subset D$ for $i \in I$ and $j \in \{1,\dots,k_i\}$, and since $D$ is convex, 
\begin{equation}\label{eq_cup_gpcs_1}
{\rm conv}\{u_{i,j}  \mid  i \in I, j=1,\dots,k_i \} \subset D.
\end{equation} 
 It is clear that $0^+D_i={\rm cone}\{v_{i,1},\dots,v_{i,\ell_i}\}+X_{i,0}$ for every $i \in I$.  As $D$ is the smallest closed convex set containing  $\bigcup\limits_{i=1}^mD_i$, we have ${\rm cone}\{v_{i,j}  \mid i \in I,  j=1,\dots,\ell_i \} \subset 0^+D$ and $X_{1,0}+\dots+X_{m,0} \subset 0^+D$. Since the cone $0^+D$ is closed, $X_0 \subset 0^+D$. Thus
  \begin{equation}\label{eq_cup_gpcs_2}
{\rm cone}\{v_{i,j}  \mid i \in I,  j=1,\dots,\ell_i \} + X_{0} \subset 0^+D.
 \end{equation} 
 Combining \eqref{eq_cup_gcps}, \eqref{eq_cup_gpcs_1} with \eqref{eq_cup_gpcs_2} yields $D'  \subset D$.  Thus we have proved that $D'=D$. Since $D'$ is a gpcs, $D$ is also a gpcs.

 Now, suppose that at least one of the set $D_1, \dots, D_m$ is polyhedral convex. Then, by Lemma \ref{rep_pcs_lemma}, in the representation \eqref{eq_rep_D_i} for $D_1, \dots, D_m$ we may assume that at least one of the sets $X_{1,0},\dots, X_{m,0}$ is finite-codimensional. According to Lemma~\ref{closed_sum_two_subspaces}, $X_{1,0}+\dots+X_{m,0}$ is a closed linear subspace of finite codimension in $X$; hence ${\rm codim}X_0 < \infty$. Due to \eqref{eq_cup_gcps}, $D'$ is a pcs by  Lemma~\ref{rep_pcs_lemma}. Since $D=D'$, we see that $D$ is a pcs. 
$\hfill\Box$

\medskip
 From Theorem~\ref{cup_gcps} we obtain the following corollary. 
 \begin{Corollary}\label{union_of_gpcs}
	If a convex subset $D \subset X$ is the union of a finite number of generalized polyhedral convex sets (resp., of polyhedral convex sets) in $X$, then $D$ is generalized polyhedral convex (resp., polyhedral convex).    
\end{Corollary}

The reader is referred to \cite[Lemma~2.50]{Rockafellar_Wets_1998} for a different proof of Corollary \ref{union_of_gpcs} in the case where $X=\mathbb{R}^n$.

\medskip
It turns out that the closure of the cone generated by a gpcs is a  generalized polyhedral convex cone. Hence, the next proposition extends \cite[Theorem~19.7]{Rockafellar_1970} to a lcHtvs setting.

\begin{Proposition}\label{generated_cone}
	If a nonempty subset $D \subset X$ is generalized polyhedral convex, then $\overline{{\rm cone}D}$ is a generalized polyhedral convex cone. In addition, if $0 \in D$ then ${\rm cone}D$ is  a generalized polyhedral convex cone; hence ${\rm cone}D$ is closed.
\end{Proposition}
\noindent
{\bf Proof.} Suppose that $D$ is of the form \eqref{rep_D}. According to \cite[Theorem~2.10]{Luan_Yen_2015}, 
\begin{equation}\label{eq_coneD_1}
C:={\rm cone}\{u_i, v_j \mid i=1,\dots,k, \, j=1,\dots,\ell\}+X_0
\end{equation}
is a generalized polyhedral convex cone, that is closed. Since $C$ contains $D$, we must have $\overline{{\rm cone}D} \subset C$. From \eqref{rep_D} we see that $0^+D={\rm cone}\{v_j \mid j=1,\dots,\ell\}+X_0$ and $u_i  \in \overline{{\rm cone}D}$ for all $i=1,\dots,k$. As $\overline{{\rm cone}D}$ is a closed convex cone, from \eqref{eq_coneD_1} it follows that $C \subset \overline{{\rm cone}D}$. Thus we have shown that $\overline{{\rm cone}D}=C$. In particular, $\overline{{\rm cone}D}$ is a generalized polyhedral convex cone. 

Now, suppose that $0 \in D$. To get the equality ${\rm cone}D=C$ with $C$ being given by \eqref{eq_coneD_1}, we first observe that ${\rm cone}D\subset C$, because $C=\overline{{\rm cone}D}$. To verify that ${\rm cone}D\supset C$, take any $x\in C$. According to \eqref{eq_coneD_1}, one can find nonnegative numbers $\lambda_1,\dots, \lambda_k,$ $\mu_1,\dots, \mu_\ell$, and a vector $x_0 \in X_0$, such that $x=\sum\limits_{i=1}^k \lambda_i u_i + \sum\limits_{j=1}^{\ell} \mu_j v_j +x_0$. If $\lambda:=\sum\limits_{i=1}^k \lambda_i$ is positive, then  $\frac{1}{\lambda}x$ belongs to $ D$; so $x \in {\rm cone}D$. If $\lambda=0$, then $\lambda_1=\dots=\lambda_k=0$ and  $x=\sum\limits_{j=1}^{\ell} \mu_j v_j +x_0$; hence $x \in 0^+D$. Since $0 \in D$, this implies that $0+x$ is contained in $D$. The inclusion ${\rm cone}D\supset C$ has been proved. So we have ${\rm cone}D=C$. In particular,  ${\rm cone}D$ is a generalized polyhedral convex cone.
$\hfill\Box$

\medskip
An analogue of Proposition~\ref{generated_cone} for polyhedral convex sets can be formulated as follows.

\begin{Proposition}\label{generated_cone_pc}
	If a nonempty subset $D \subset X$ is polyhedral convex, then $\overline{{\rm cone}D}$ is a polyhedral convex cone. In addition, if $0 \in D$ then ${\rm cone}D$ is  a polyhedral convex cone; hence ${\rm cone}D$ is closed. 
\end{Proposition}
\noindent
{\bf Proof.}  Similar to the proof of Proposition \ref{generated_cone}. $\hfill\Box$

\medskip
In convex analysis, to every convex set and a point belonging to it, one associates a tangent cone. Let us complete this section by showing that the tangent cone to a gpcs at a given point is a generalized polyhedral convex cone. By definition, the (Bouligand-Severi) \textit{tangent cone} \cite{Aubin_Frankowska_1990} $T_D(x)$ to a closed subset $D \subset X$ at $x \in D$ is the set of all $v\in X$ such that there exist sequences $t_k\to 0^+$ and $v_k\to v$ such that $x+t_kv_k\in D$ for every $k$. If $D$ is convex, then 
\begin{equation}\label{eq_tangent_cone}
T_D(x)=\overline{{\rm cone}(D-x)}.
\end{equation} 
If $D$ is a gpcs and $x \in D$, then $D-x$ is a gpcs containing $0$. Therefore, according to Proposition~\ref{generated_cone}, ${\rm cone}(D-x)$ is a generalized polyhedral convex cone, that is closed. So the closure sign in the right-hand side of \eqref{eq_tangent_cone} can be omitted.  Similarly, according to Proposition~\ref{generated_cone_pc}, if $D$ is a pcs and $x \in D$, then ${\rm cone}(D-x)$ is a polyhedral convex cone and the closure sign in the right-hand side of \eqref{eq_tangent_cone} can be also omitted. Thus we have obtained the following result. 
\begin{Proposition}\label{tangent_cone}
	If $D \subset X$ is a generalized polyhedral convex set (resp., a polyhedral convex set) and if  $x \in D$, then $T_D(x)$ is a generalized polyhedral convex cone (resp., a polyhedral convex cone) and one has 
	$T_D(x)={\rm cone}(D-x)$.
\end{Proposition}

\section{Generalized Polyhedral Convex Functions}
\setcounter{equation}{0}  

As the title indicates, this section will deal with the concept of generalized polyhedral convex function. The latter is based on the notion of generalized polyhedral convex set which has been considered in details in the preceding section.

\medskip
Let $f$ be a function from a locally convex Hausdorff topological vector space $X$ to $\overline{\mathbb{R}}:=\mathbb{R} \cup \{\pm \infty\}$. The \textit{effective domain} and the \textit{epigraph} of $f$ are defined, respectively, by setting ${\rm dom}f=\{x \in X \mid f(x) < +\infty\}$ and $${\rm epi}f=\big\{(x, \alpha) \in X \times \mathbb{R} \mid x \in {\rm dom}f,\ f(x) \leq \alpha \big\}.$$  If ${\rm dom}f$ is nonempty and $f(x) > - \infty$ for all $x \in X$, then $f$ is said to be \textit{proper}.  We say that $f$ is  \textit{convex} if  ${\rm epi}f$ is a convex set in $X \times \mathbb{R}$. It is easily verified that $f$ is convex if and only if ${\rm dom}f$ is convex and the Jensen inequality
$$f((1-t)x_1+tx_2)\leq (1-t)f(x_1)+tf(x_2)$$ is valid for any $x_1,\, x_2$ in ${\rm dom}f$ and $t\in (0,1)$.

\medskip
According to Rockafellar \cite[p.~172]{Rockafellar_1970}, a real-valued function defined on $\mathbb{R}^n$ is called polyhedral convex if its epigraph is a polyhedral convex set in $\mathbb{R}^{n+1}$. The following notion of generalized polyhedral convex function appears naturally in that spirit.

\begin{Definition} {\rm 
	 Let $X$ be a locally convex Hausdorff topological vector space. A function $f:X\to\overline{\mathbb{R}}$ is called \textit{generalized polyhedral convex} (resp., \textit{polyhedral convex}) if its epigraph is a generalized polyhedral convex set (resp., a polyhedral convex set) in $X\times\mathbb R$.}
\end{Definition}

Complete characterizations of a generalized polyhedral convex function (resp., a polyhedral convex function) in the form of the maximum of a finite family of continuous affine functions over a certain generalized polyhedral convex set (resp., a polyhedral convex set) are given in the next theorem.

\begin{Theorem}\label{rep_gpcf}	Suppose that $f:X\to\overline{\mathbb{R}}$ is a proper function. Then $f$ is generalized polyhedral convex (resp., polyhedral convex) if and only if ${\rm dom}f$ is a generalized polyhedral convex set (resp., a polyhedral convex set) in $X$ and there exist $v_k^* \in X^*$, $\beta_k \in \mathbb{R}$, for $k=1,\dots,m$, such that
		\begin{equation}\label{eq_rep_gcpf}
		f(x)=\begin{cases}
		\max \big\{ \langle v_k^*, x \rangle + \beta_k \mid k=1,\dots,m  \big\} &\text{if } x \in {\rm dom}f,\\
		+\infty & \text{if } x \notin {\rm dom}f.
		\end{cases}
		\end{equation}	  
\end{Theorem}
\noindent
{\bf Proof.} Let $f:X\to\overline{\mathbb{R}}$ be a proper function.

First, suppose that $f$ is a gpcf. Then there exist a closed affine subspace $L \subset X\times \mathbb{R}$,\, $u^*_i \in X^*,\, a_i \in \mathbb{R},\, b_i \in \mathbb{R}$, for $i=1,\dots, p,$ such that
\begin{equation}\label{proof_epi_f}
{\rm epi}f=\left\{(x, t) \in L \mid \langle u^*_i, x \rangle +a_i t\leq b_i,\  i=1,\dots,p\right\}.
\end{equation}
By \cite[Remark 2.196]{Bonnans_Shapiro_2000}, one can find a continuous linear mapping $\widetilde{A}$ from $X\times \mathbb{R}$ to a lcHtvs $Y$ and $y \in Y$ so that $L=\big\{(x, t) \in X\times \mathbb{R} \mid \widetilde{A}(x,t)=y  \big\}$. Let the continuous linear mapping $A: X \to Y$ be defined by $A(x)= \widetilde{A}(x,0)$. For $y_0:=\widetilde{A}(0,1)$, we see that
\begin{equation}\label{tildeA(x,t)}
\widetilde{A}(x,t) = \widetilde{A}(x,0)+t\widetilde{A}(0,1)=A(x)+ty_0\quad  (x \in X,\, t \in \mathbb{R}).
\end{equation}
Given any $(\bar{x}, \bar{t}) \in {\rm epi}f$, since $(\bar{x}, \bar{t}+\gamma) \in {\rm epi}f$ for all $\gamma \geq 0$. In particular,  $(\bar{x}, \bar{t}+\gamma) \in L$ for all $\gamma \geq 0$. So we have 
\begin{equation*}
  y= A\bar{x}+(\bar{t}+\gamma) y_0 = (A\bar{x}+ \bar{t} y_0) + \gamma y_0= y+ \gamma y_0    
\end{equation*}
for all $\gamma \geq 0$. It follows that $y_0=0$. Substituting $(x,t)=(\bar{x}, \bar{t}+\gamma)$ into the inequalities in \eqref{proof_epi_f} yields 
 $a_i \leq 0$ for $i=1,\dots,p$. There exists an index $i \in \{1,\dots,p\}$ satisfying $a_i <0$. Indeed, suppose on the contrary that $a_i =0$  for $i=1,\dots,p$. By the properness of $f$, there exists $\bar x\in X$ with $|f(\bar x)|<\infty$. Then $(\bar x, f(\bar x))\in {\rm epi}f\subset L$. As $y_0=0$, from \eqref{tildeA(x,t)} it follows that $y=\widetilde{A}(\bar x,f(\bar x))=A(\bar x)$. Moreover, for any $t\in\mathbb R$, combining this with \eqref{tildeA(x,t)} one has $\widetilde{A}(\bar x,t)=y$. Hence $(\bar x,t)\in L$ for all $t\in\mathbb R$. Since $\langle u^*_i,\bar x \rangle +a_i t=\langle u^*_i, \bar x \rangle +a_i f(\bar x)\leq b_i$ for $i=1,\dots,p$, we see that  $(\bar x,t)\in {\rm epi}f$ for all~$t\in\mathbb R$. Then $f(\bar x) =-\infty$. We have thus arrived at a contradiction.
 
For each $i \in \{1,\dots,p\}$ with $a_i<0$, we replace the inequality $\langle u^*_i, x \rangle +a_i t\leq b_i$  by the following equivalent one: 
$$\left\langle\frac{1}{|a_i|}u^*_i, x \right\rangle - t\leq \frac{b_i}{|a_i|}.$$
Then, reordering the family $\{a_1,\dots,a_p\}$ (if necessary), we may assume that $a_k=-1$ for $k=1,\dots,m,$ with $m \leq p$, and $a_i=0$ for $i=m+1,\dots,p$. It follows that
 \begin{equation}\label{rep_epi_f}
 \begin{aligned}
 {\rm epi}f=\big\{(x, t) \in X \times \mathbb{R} \mid Ax=y, \ & \langle u^*_k, x \rangle -b_k \leq t,\, k=1,\dots,m,\\
 &\; \; \; \langle u^*_i, x \rangle \leq b_i,\,  i=m+1,\dots,p\big\}.
 \end{aligned}
 \end{equation} 
This implies that
\begin{equation}\label{rep_dom_f}
{\rm dom}f=\left\{x \in X \mid Ax=y, \langle u^*_i, x \rangle \leq b_i, \, i=m+1,\dots,p\right\}.
\end{equation}
In particular, ${\rm dom}f$ is a gpcs in $X$. Combining \eqref{rep_epi_f} with \eqref{rep_dom_f} gives 
\begin{equation*}\label{2nd_rep_epi_f}
{\rm epi}f=\big\{(x, t) \in X \times \mathbb{R} \mid x \in  {\rm dom}f,\; \langle u^*_k, x \rangle -b_k \leq t,\, k=1,\dots,m\big\}.
\end{equation*}
So, $f$ can be represented in the form \eqref{eq_rep_gcpf} with $v^*_k:=u^*_k$, $\beta_k:=-b_k$ for $k=1,\dots,m$. 

In addition, if $f$ is a pcf on $X$, then we may assume that ${\rm epi}f$ is of the form \eqref{proof_epi_f}, where $L=X \times \mathbb{R}$. In this case, we can repeat the above proof with $Y:=\{0\}$ (the trivial space), $\widetilde{A}(x,t)\equiv 0$ and $y=0$. Hence, it follows from \eqref{rep_epi_f} that $f$ admits the representation \eqref{eq_rep_gcpf} with ${\rm dom}f$ being a pcs in $X$. 

Now, suppose that ${\rm dom}f$ is a gpcs in $X$ and $f$ is given by \eqref{eq_rep_gcpf}. Then there exist $x^*_i \in X^*$, $\alpha_i \in \mathbb{R}$, $i=1,\dots,p,$ a continuous linear mapping $B$ from $X$ to a lcHtvs $Z$, and a vector $z \in Z$ such that
\begin{equation}\label{rep_dom_f_2}
{\rm dom}f=\big\{x \in X \mid Bx=z,\; \langle x^*_i, x \rangle \leq \alpha_i,\; i=1,\dots,p\big \}.
\end{equation}
Combining this with \eqref{eq_rep_gcpf}, we obtain
\begin{equation}\label{rep_epi_f_2}
	\begin{aligned}
	{\rm epi}f	&=\big\{(x, t) \in X \times \mathbb{R} \mid Bx=z, \, \langle x^*_i, x \rangle \leq \alpha_i, \, i=1,\dots,p, \\
		&\hspace*{5.4cm}\langle v^*_k, x \rangle + \beta_k \leq t, \ k=1,\dots,m \big\}\\
		&=\big\{(x, t) \in X \times \mathbb{R} \mid Bx+0t=z, \langle x^*_i, x \rangle +0t\leq \alpha_i,\, i=1,\dots,p, \\
		&\hspace*{6.1cm}
	\langle v^*_k, x \rangle - t \leq -\beta_k ,\, k=1,\dots,m \big\}.\\
	\end{aligned}
\end{equation}     
This clearly shows that ${\rm epi}f$ is a gpcs in $X \times \mathbb{R}$; hence $f$ is a gpcf. 

Finally, let us assume that ${\rm dom}f$ is a pcs in $X$ and $f$ is given by \eqref{eq_rep_gcpf}. Then, in the formula \eqref{rep_dom_f_2} for ${\rm dom}f$, we can choose $Z=\{0\}$, $B\equiv 0$, and $z=0$. Clearly, with the chosen $B$ and $z$, \eqref{rep_epi_f_2} implies that  ${\rm epi}f$ is a pcs; so $f$ is polyhedral convex. 
$\hfill\Box$

\begin{Remark}
{\rm For the case $X=\mathbb{R}^n$, the result in Theorem \ref{rep_gpcf} is a known one (see \cite[p.~172]{Rockafellar_1970}, \cite[Theorem~2.49]{Rockafellar_Wets_1998}, \cite[Proposition~3.2.3]{Bertsekas_et_al_2003}). In the above proof, we have used some ideas of the proof of \cite[Proposition~3.2.3]{Bertsekas_et_al_2003}.}
\end{Remark}

Theorem \ref{rep_gpcf} provides us with a general formula for any generalized polyhedral convex function on a locally convex Hausdorff topological vector space. For polyhedral convex functions on $\mathbb{R}^n$, there is another important characterization \cite[Theorem~2.49]{Rockafellar_Wets_1998}: \textit{A proper convex function $f$ is polyhedral convex if and only if $f$ is piecewise linear}. In order to obtain an analogous result for generalized polyhedral convex functions, we need the following infinite-dimensional generalization of the concept of piecewise linear function on $\mathbb{R}^n$ of \cite{Rockafellar_Wets_1998}.

\begin{Definition}\label{def_gpl} {\rm A proper function $f:X \to \overline{\mathbb{R}}$, which is defined on a locally convex Hausdorff topological vector space, is said to be \textit{generalized piecewise linear} (resp., \textit{piecewise linear}) if there exist generalized polyhedral convex sets (resp., polyhedral convex sets) $D_1, \dots, D_m$  in $X$, $v_1^*, \dots, v_m^* \in X^*$, and $\beta_1, \dots, \beta_m \in \mathbb{R}$ such that  ${\rm dom}f=\bigcup\limits_{k=1}^m D_k$ and 
$f(x)=\langle v_k^*, x \rangle + \beta_k$ for all $x \in D_k$, $k=1,\dots,m$.}\end{Definition}

\begin{Remark}\label{sum_of_plfs}
{\rm By using the definition, it is not difficult to show that the sum a finite family of generalized piecewise linear functions (resp., a finite family of piecewise linear functions) is a generalized piecewise linear function (resp., piecewise linear function).}
\end{Remark}

The forthcoming theorem clarifies the relationships between generalized polyhedral convex functions and generalized piecewise linear functions. 

\begin{Theorem}\label{rep_gpcf_gpl}
	A proper convex function is generalized polyhedral convex (resp., polyhedral convex) if and only if it is generalized piecewise linear (resp., piecewise linear). 
\end{Theorem}
\noindent
{\bf Proof.} Let $f:X\to\overline{\mathbb{R}}$ be a proper convex function.

First, suppose that $f$ is generalized polyhedral convex (resp., polyhedral convex). By Theorem~\ref{rep_gpcf}, ${\rm dom}f$ is a gpcs (resp., a pcs) and there exist $v_k^* \in X^*$, $\beta_k \in \mathbb{R}$, $k=1,\dots,m$, such that $f(x)=\max \left\{ \langle v_k^* , x \rangle + \beta_k  \mid k=1,\dots,m  \right\}$  for all $x \in {\rm dom}f$.
For each $k=1,\dots,m$, put  
\begin{equation*}
\begin{aligned}
D_k &={\rm dom}f  \cap \big\{ x \in X \mid  \langle v_i^* , x \rangle + \beta_i \leq \langle v_k^* , x \rangle + \beta_k,\; \forall i=1,\dots,m  \big\}\\
& ={\rm dom}f  \cap \big\{ x \in X \mid  \langle v_i^*- v_k^*, x \rangle \leq \beta_k-\beta_i,\; \forall i=1,\dots,m  \big\}.
\end{aligned}
\end{equation*}
Observe that ${\rm dom}f=\bigcup\limits_{k=1}^m D_k$ and $f(x)= \langle v_k^* , x \rangle + \beta_k$ for every $x \in D_k$, $k=1,\dots,m$. Since ${\rm dom}f$ is a gpcs (resp., a pcs), $D_k$ is also a gpcs (resp., a pcs). It follows that $f$ is a generalized piecewise linear function (resp., a piecewise linear function).  

Now, suppose that $f$ is generalized piecewise linear (resp., piecewise linear). Then, one can find  generalized polyhedral convex sets (resp., polyhedral convex sets) $D_1, \dots, D_m$ in $X$, $v_1^*, \dots, v_m^* \in X^*$, and $\beta_1, \dots, \beta_m \in \mathbb{R}$ such that  ${\rm dom}f=\bigcup\limits_{k=1}^m D_k$ and $f(x)=\langle v^*_k, x \rangle + \beta_k$ for all $x \in D_k$, $k=1,\dots,m$. It follows that ${\rm epi}f=\bigcup\limits_{k=1}^m E_k$, where
$$E_k:=\big\{(x, t) \in X \times \mathbb{R} \mid x \in D_k, \, \langle v_k^*, x \rangle + \beta_k\leq t \big \} \quad (k=1,\dots,m).$$
So, for each $k =1,\dots,m$, $E_k$ is the intersection of the generalized polyhedral convex set (resp.,  the polyhedral convex set) $D_k \times \mathbb{R}$ and the polyhedral convex set $$\big\{(x, t) \in X \times \mathbb{R} \mid \langle v_k^*, x \rangle + \beta_k\leq t \big \}.$$ In particular, $E_k$ is a gpcs (resp., a pcs). The convexity of $f$ shows that ${\rm epi}f$ is convex. Combining this with the fact that ${\rm epi}f$ is the union of the gpcs (resp., the pcs) $E_1, \dots, E_m$, we conclude by Corollary~\ref{union_of_gpcs} that the set ${\rm epi}f$ is generalized polyhedral convex (resp.,  polyhedral convex). Thus $f$ is a generalized polyhedral convex function (resp., a polyhedral convex function). 

The proof is complete. 
$\hfill\Box$

\medskip
Based on Theorem~\ref{rep_gpcf_gpl}, we can prove that the class of generalized polyhedral convex functions (resp., the class of polyhedral convex functions) is invariant w.r.t. the addition of functions. 

\begin{Theorem}\label{sum_two_gcpf} 
	Let $f_1, f_2$ be two proper functions on $X$. If $f_1, f_2$ are generalized polyhedral convex (resp., polyhedral convex) and $({\rm dom}f_1)\cap ({\rm dom}f_2)$ is nonempty, then $f_1+f_2$ is a proper generalized polyhedral convex function (resp., a polyhedral convex function).    
\end{Theorem}
\noindent
{\bf Proof.} Suppose that $f_1, f_2$ are proper generalized polyhedral convex functions (resp., proper polyhedral convex functions) defined on $X$ with  $({\rm dom}f_1)\cap ({\rm dom}f_2)\neq\emptyset$. Then, $f_1+f_2$ is a proper convex function. Due to Theorem~\ref{rep_gpcf_gpl}, $f_1$ and $f_2$ are generalized piecewise linear (resp., piecewise linear); hence $f_1+f_2$ is generalized piecewise linear (resp., piecewise linear) by Remark~\ref{sum_of_plfs}. So, according to Theorem~\ref{rep_gpcf_gpl}, the function $f_1+f_2$ is generalized polyhedral convex (resp., polyhedral convex). 
$\hfill\Box$

\begin{Remark}
	{\rm For the case $X=\mathbb{R}^n$, the result in Theorem~\ref{sum_two_gcpf} is a known one (see \cite[Theorem~19.4]{Rockafellar_1970}).}
\end{Remark}

In convex analysis, it is well known \cite{Ioffe_Tihomirov_1979, Rockafellar_1970} that the concept of directional derivative has an important role. We are going to discuss a property of the directional derivative mapping of a gpcf (resp., a pcf) at a given point. 

\medskip
If $f: X \to \overline{\mathbb{R}}$ is a proper convex function and $x \in X$ which $f(x)$ is finite, the \textit{directional derivative} of $f$ at $x$ w.r.t. a direction $h \in X$, denoted by $f'(x; h)$, always exists (it can take values $-\infty$ or $+\infty$). For the case where $X=\mathbb{R}^n$, the proof can be found in \cite[Theorem~23.1]{Rockafellar_1970}. For the case where $X$ is a lcHtvs, the fact has been discussed, e.g., in \cite[pp.~48--49]{Bonnans_Shapiro_2000} and \cite[Theorem~2.1.13]{Zalinescu_2002}. According to \cite[Proposition~2.60]{Bonnans_Shapiro_2000}, the closure of the epigraph of $f'(x; \cdot)$ coincides with the tangent cone to ${\rm epi}f$ at $(x, f(x))$, i.e.,
\begin{equation}\label{epi_dir_der_1}
\overline{{\rm epi}f'(x; \cdot)}=T_{{\rm epi}f}(x, f(x)).
\end{equation}

By \cite[Theorem~23.10]{Rockafellar_1970} we know that if $f: \mathbb{R}^n \to \overline{\mathbb{R}}$ is proper polyhedral convex, then the closure sign in \eqref{epi_dir_der_1} can be omitted and $f'(x; \cdot)$ is a proper polyhedral convex function. The last two facts can be extended to polyhedral convex functions on lcHtvs and generalized polyhedral convex functions as follows.

\begin{Theorem}\label{directional_derivative}
	Let $f$ be a proper generalized polyhedral convex function (resp., a proper polyhedral convex function) on a locally convex Hausdorff topological vector space $X$. For any $x \in {\rm dom}\,f$, $f'(x; \cdot)$ is a proper generalized polyhedral convex function (resp., a proper polyhedral convex function). In particular, ${\rm epi}f'(x; \cdot)$ is closed and, by \eqref{epi_dir_der_1} one has 
	\begin{equation}\label{epi_dir_der_2}
	{\rm epi}f'(x; \cdot)=T_{{\rm epi}f}(x, f(x)).
	\end{equation}
\end{Theorem}

For proving this theorem, we need the next lemma, which can be proved by using definition of tangent cone. For the case $X$ is a Banach space, formula \eqref{eq_tangent_cone_of_sum} has been given in \cite[Proposition~3.1]{Ban_Mordukhovich_Song_2011}.

\begin{Lemma}\label{tangent_cone_of_sum}
	Suppose that a subset $D \subset X$ is the union of nonempty closed convex sets $D_1,\dots,D_m$. Then, for every $x \in D$, 
	\begin{equation*}\label{eq_tangent_cone_of_sum}
	T_D(x)=\bigcup\limits_{k \in J(x)} T_{D_k}(x),
	\end{equation*}
	where $J(x)=\big\{j \in \{1,\dots,m\} \mid x \in D_j\big\}$.  
\end{Lemma}

\noindent
{\bf Proof of Theorem~\ref{directional_derivative}.} Let $f$ be given as in the formulation of the theorem. Due to the ``only if" part of Theorem~\ref{rep_gpcf_gpl}, one can find  generalized polyhedral convex sets (resp., polyhedral convex sets) $D_1, \dots, D_m$ in $X$, $v_1^*, \dots, v_m^* \in X^*$, and $\beta_1, \dots, \beta_m \in \mathbb{R}$ such that  ${\rm dom}f=\bigcup\limits_{k=1}^m D_k$ and $f(x)=\langle v^*_k, x \rangle + \beta_k$ for all $x \in D_k$, $k=1,\dots,m$. We may assume that $D_k\neq\emptyset$ for all $k\in \{1,\dots,m\}$. Given any $x \in {\rm dom}f$, $f(x)$ is finite because $f$ is proper. By \cite[p.~195]{Ioffe_Tihomirov_1979}, $f'(x; \cdot)$ is convex. 

If $h\notin T_{{\rm dom}f}(x)$, then $x+th \notin {\rm dom}f$ for every $t>0$; so $f'(x,h)=+\infty$. If $h \in T_{{\rm dom}f}(x)$, then by Lemma~\ref{tangent_cone_of_sum} one can find $D_k$ such that $x \in D_k$ and $h \in T_{D_k} (x)$. In addition, since $D_k$ is a gpcs (resp., a pcs), $T_{D_k}(x)={\rm cone}(D_k-x)$ by Proposition~\ref{tangent_cone}. Hence, there exists $\delta>0$ satisfying $x+\delta h \in D_k$. As $D_k$ is convex, one has  $x+th \in D_k$ for all $t \in [0, \delta]$. It follows that
\begin{equation*}
\frac{f(x+th)-f(x)}{t}=\frac{\big(\langle v^*_k, x+th \rangle + \beta_k \big) -\big(\langle v^*_k, x \rangle + \beta_k \big)}{t}=\langle v^*_k,h \rangle
\end{equation*}
for every $t \in (0, \delta]$; so $f'(x;h)=\langle v^*_k,h \rangle$. Note that, if $h \in T_{D_{k_1}} (x) \cap T_{D_{k_2}} (x)$, then $\langle v^*_{k_1},h \rangle=\langle v^*_{k_2},h \rangle.$ Indeed, one can find positive numbers $\delta_1, \delta_2$ such that $x+th \in D_{k_1}$ for all $t \in [0, \delta_1]$ and $x+th \in D_{k_2}$ for all $t \in [0, \delta_2]$. Setting $\delta=\min\{\delta_1, \delta_2\}$, we have $x+th \in D_{k_1} \cap D_{k_2}$ for every $t \in [0, \delta]$. Then
\begin{equation}\label{proof_dir_der}
f(x+th)=\langle v_{k_1}^*, x +th \rangle + \beta_{k_1}=\langle v_{k_2}^*, x +th \rangle + \beta_{k_2} \quad (t \in [0, \delta]).
\end{equation}  
 In particular, $\langle v_{k_1}^*, x\rangle + \beta_{k_1}=\langle v_{k_2}^*, x\rangle + \beta_{k_2}$. Hence, \eqref{proof_dir_der} implies $\langle v^*_{k_1},h \rangle=\langle v^*_{k_2},h \rangle.$ 
 
 We have shown that $f'(x; h)$ is well defined and finite for every $h \in T_{{\rm dom}f}(x)$. Consequently, the function $f'(x;\cdot)$ is proper and ${\rm dom}f'(x; \cdot)=T_{{\rm dom}f}(x)$. Moreover, applying Lemma~\ref{tangent_cone_of_sum} for $D={\rm dom}f$, we can assert that $T_{{\rm dom}f}(x)$ is the union of the generalized polyhedral convex cones (resp., the union of the polyhedral convex cones) $T_{D_k}(x)$, ${k \in J(x)}$, and $f'(x; h)= \langle v^*_{k},h \rangle$ if $h \in T_{D_k}(x)$ with ${k \in J(x)}$. This implies that the proper convex function $f'(x; \cdot)$ is generalized piecewise linear (resp., piecewise linear). Hence, by the ``if" part of Theorem~\ref{rep_gpcf_gpl}, $f'(x; \cdot)$ is proper generalized polyhderal convex (resp., proper polyhderal convex). Therefore, ${\rm epi}f'(x; \cdot)$ is closed and \eqref{epi_dir_der_2} follows from \eqref{epi_dir_der_1}. $\hfill\Box$

\medskip
In the final part of this section, we are interested in the concept of infimal convolution function, which was introduced by Fenchel \cite{Fenchel_1953} and discussed by many other authors (see, e.g., Rockafellar \cite{Rockafellar_1970}, Ioffe and Tihomirov \cite{Ioffe_Tihomirov_1979}, Attouch and Wets \cite{Attouch_Wets_1989}, Str\"{o}mberg \cite{Stromberg_1994, Stromberg_1996}). 
According to Rockafellar \cite[p.~34]{Rockafellar_1970}, the infimal convolution operation is analogous to the classical formula for integral convolution and, in a sense, is dual to the operation of addition of convex functions. 

\medskip
As noted by Nam \cite[p.~2215]{Nam_2015} and Nam and Cuong \cite[pp.~333--334]{Nam_Cuong_2015}, a large spectrum of known nonsmooth functions can be interpreted as infimal convolutions. In the above cited papers, the authors have obtained some upper estimates for three types of subdifferentials of a class of nonconvex infimal convolutions. 

\medskip
Although the infimal convolution of a finite family of functions can be defined \cite{Ioffe_Tihomirov_1979, Rockafellar_1970}, for simplicity, we will only consider the infimal convolution of two functions.

\begin{Definition}\label{def_infimal_convol} {\rm (See   \cite[p.~168]{Ioffe_Tihomirov_1979} and \cite[p.~34]{Rockafellar_1970})
	Let $f_1, f_2$ be two proper functions on a locally convex Hausdorff topological vector space $X$. The \textit{infimal convolution} of $f_1, f_2$ is the function defined by
	\begin{equation}\label{def_conv_func}
	(f_1 \square f_2)(x):=\inf\left\{f_1(x_1)+f_2(x_2) \mid x_1+x_2=x\right\}.
	\end{equation}  
}\end{Definition}

If $f_1, f_2$ are proper convex, then $f_1 \square f_2$ is convex (see, e.g., \cite[p.~43]{Zalinescu_2002}). However, if $f_1, f_2$ are proper, $f_1 \square f_2$ may not be proper. For instance, choosing $X=\mathbb{R}$, $f_1(x)=x$ and $f_2(x)=2x$, one has $(f_1 \square f_2)(x)=-\infty$ for all $x \in X$. 

\medskip
The infimal convolution operation in \eqref{def_conv_func} corresponds to the addition of the epigraphs of $f_1$ and $f_2$  as sets in $X \times \mathbb{R}$. Namely, as noted in \cite[p.~34]{Rockafellar_1970},
\begin{equation*}\label{def_conv_func_1}
(f_1 \square f_2)(x)=\inf\left\{\alpha \mid (x,\alpha)\in {\rm epi}f_1+{\rm epi}f_2\right\}.
\end{equation*}   

 According to Proposition~\ref{sum_pc_gpc}, the sum of a polyhedral convex set and a generalized polyhedral convex set is a polyhedral convex set. We will use this fact to prove the following proposition.

\begin{Proposition}\label{inf_convol_pcs_gpcs} Let $f_1, f_2$ be two proper functions. If $f_1$ is polyhedral convex and $f_2$ is generalized polyhedral convex, then $f_1 \square f_2$ is a polyhedral convex function.  
\end{Proposition}
\noindent
{\bf Proof.} First, let us verify the inclusion ${\rm epi}f_1+{\rm epi}f_2 \subset {\rm epi}(f_1 \square f_2)$. Pick any $(x_i, \alpha_i) \in {\rm epi}f_i$, $i=1,2$. Then we have $f_1(x_1)+f_2(x_2) \leq \alpha_1+\alpha_2$. Combining this with \eqref{def_conv_func}, we get $(f_1 \square f_2)(x_1+x_2) \leq \alpha_1+\alpha_2$; hence $(x_1+x_2, \alpha_1+\alpha_2) \in {\rm epi}(f_1 \square f_2)$. 

Now, to show that ${\rm epi}(f_1 \square f_2) \subset {\rm epi}f_1+{\rm epi}f_2$, select a point $(x, \alpha) \in {\rm epi}(f_1 \square f_2)$. For any $\varepsilon >0$, since $(f_1 \square f_2)(x) \leq \alpha $, there exist $x_1 \in {\rm dom}f_1$  and $x_2 \in {\rm dom}f_2$ such that $x_1+x_2=x$ and $f_1(x_1)+f_2(x_2) \leq \alpha+\varepsilon$. As $f_2(x_2) \leq \alpha+\varepsilon - f_1(x_1)$, one has
$$(x, \alpha+\varepsilon)=\left(x_1, f_1(x_1)\right)+\left(x_2, \alpha+\varepsilon - f_1(x_1)\right) \in {\rm epi}f_1+{\rm epi}f_2.$$  
So, letting $\varepsilon \to 0^+$ yields $(x, \alpha) \in \overline{{\rm epi}f_1+{\rm epi}f_2}$. Since ${\rm epi}f_1$ is a pcs and ${\rm epi}f_2$ is a gpcs, ${\rm epi}f_1+{\rm epi}f_2$ is a pcs by Proposition~\ref{sum_pc_gpc}. In particular, ${\rm epi}f_1+{\rm epi}f_2$ is closed. Hence, $(x, \alpha) \in {\rm epi}f_1+{\rm epi}f_2$. 

We have thus proved that ${\rm epi}(f_1 \square f_2)= {\rm epi}f_1+{\rm epi}f_2$, where the set on the right-hand side is polyhedral convex. This means that $f_1 \square f_2$ is a polyhedral convex function.  $\hfill\Box$

\begin{Remark}{\rm Proposition~\ref{inf_convol_pcs_gpcs} is a generalization of \cite[Corollary~19.3.4]{Rockafellar_1970}), where the case $X=\mathbb{R}^n$ was treated. If $X$ is a general locally convex Hausdorff topological vector space and $f_1, f_2$ are generalized polyhedral convex, $f_1 \square f_2$ may not be generalized polyhedral convex. To see this, one can choose a suitable space $X$ and closed linear subspaces $X_1, X_2$ of $X$ such that  $\overline{X_1+X_2}=X$ and $X_1+X_2 \neq X$ (see Remark \ref{Rem_sum_closed_spaces} for details). Let $f_i:=\delta(\cdot, X_i)$ ($i=1,2$) be the \textit{indicator function} of $X_i$, i.e., $f_i(x)=0$ for $x\in X_i$ and $f_i(x)=+\infty$ for $x\notin X_i$. Clearly, both functions $f_1$ and $f_2$ are proper generalized polyhedral convex. An easy computation shows that $(f_1 \square f_2)(\cdot)= \delta(\cdot, X_1+X_2)$ and ${\rm epi}(f_1 \square f_2)=(X_1+X_2)\times [0,+\infty)$. Since $X_1+X_2$ is non-closed, ${\rm epi}(f_1 \square f_2)$ is non-closed; hence $f_1 \square f_2$ cannot be a gpcf.} 
\end{Remark}

\section{Dual Constructions}
\setcounter{equation}{0}  

Various properties of normal cones to and  polars of generalized polyhedral convex sets,  conjugates of generalized polyhedral convex functions, and subdifferentials of generalized polyhedral convex functions will be studied in this section. As before, $X$ is a locally convex Hausdorff topological vector space and $X^*$ is the dual space of $X$. According to \cite[Theorem~3.10]{Rudin_1991} (see also the property of the dual space described in \cite[p.~65]{Rudin_1991}), the weak$^*$--topology makes $X^*$  into a locally convex Hausdorff topological vector space whose dual space is $X$. 

\medskip
Now, suppose that $C \subset X$ is a nonempty convex set. The \textit{normal cone} \cite[p.~205]{Ioffe_Tihomirov_1979} to $C$ at $x \in C$ is the set $N_C(x):=\big\{x^* \in X^* \mid \langle x^*, u-x\rangle \leq 0, \, \forall u \in C\big\}.$ The formula $C^{\perp}:=\{x^* \in X^* \mid \langle x^*, u \rangle = 0, \ \forall u \in C\}$ defines the \textit{annihilator} {\rm \cite[p.~117]{Luenberger_1969}} of $C$. Note that $N_C(x)$ is a closed convex cone in $X^*$, while $C^{\perp}$ is a closed linear subspace of $X^*$. If $C$ is a linear subspace of $X$, then $N_C(x)$ does not depend on $x$. Moreover, $N_C(x)=C^{\perp}$ for all $x \in C$. 
 
\medskip
In the sequence, if not otherwise stated, $D \subset X$ is a nonempty generalized polyhedral convex set given by \eqref{eq_def_gpcs_2}. Set $I=\{1,\dots,p\}$ and $I(x)=\left\{i \in I \mid \langle x^*_i, x \rangle = \alpha_i \right\}$ for $x \in D$. If $D$ is a pcs, then one can choose $Y=\{0\}$, $A\equiv 0$, and $y=0$. 

\medskip
Normal cones to a gpcs also share the polyhedrality structure. 

\begin{Theorem}\label{normal_cone_gpcs} 
	If $D \subset X$ is a generalized polyhedral convex set and if  $x \in D$, then $N_D(x)$ is a generalized polyhedral convex cone.
\end{Theorem}
\noindent
{\bf Proof.} Since $({\rm ker}A)^{\perp} \subset X^*$ is a closed linear subspace, by using \cite[Theorem 2.10]{Luan_Yen_2015} we can assert that
\begin{equation*}
Q_x:={\rm cone}\{x^*_i \mid i \in I(x)\}+({\rm ker}A)^{\perp}
\end{equation*}
is a generalized polyhedral convex cone of $X^*$. In particular, $Q_x$ is convex and closed. To show that $Q_x \subset N_D(x)$, take any $x^* \in Q_x$. Then there exist $\lambda_i \geq 0$ for $i \in I(x)$ and $u^* \in ({\rm ker}A)^{\perp}$ such that $x^*=\sum\limits_{i \in I(x) }\lambda_i x^*_i+u^*$. For any $u\in D$, from \eqref{eq_def_gpcs_2} it follows that $\langle u^*, u-x \rangle=0$, because $u-x \in {\rm ker}A$. Hence
\begin{equation*}
\begin{aligned}
\langle x^*, u-x \rangle &= \sum\limits_{i \in I(x) }\lambda_i \big( \langle x^*_i, u \rangle -\langle x^*_i, x \rangle \big) +\langle u^*, u-x \rangle\\
&=\sum\limits_{i \in I(x) }\lambda_i \big(\langle x^*_i, u \rangle -\alpha_i \big) \leq 0.
\end{aligned}
\end{equation*}
The last inequality is clear as $\lambda_i \geq 0$ for all $i \in I(x)$ and $u \in D$. Hence $\langle x^*, u-x \rangle \leq 0$ for every $u \in D$; so $x^* \in N_D(x)$. We have thus proved that $Q_x \subset N_D(x)$. To obtain the opposite inclusion, take any $v^* \in X^* \setminus Q_x$. Since $\{v^*\} \cap Q_x=~\emptyset$, by the strong separation theorem \cite[Theorem~3.4 (b)]{Rudin_1991} we can find $v \in X$ such that \begin{equation}\label{separation_Qx}
\langle v^*, v \rangle > \sup\{\langle x^*,v \rangle \mid  x^* \in Q_x\}.
\end{equation} As the linear functional $\langle \cdot,v\rangle $ is bounded from the above on the generalized polyhedral convex set $Q_x={\rm cone}\{x^*_i \mid i \in I(x)\}+({\rm ker}A)^{\perp}$, by \cite[Theorem~3.3]{Luan_Yen_2015} we know that the linear programming problem
$\max\{\langle x^*,v \rangle \mid  x^* \in Q_x\}$ has a solution. Therefore,
 invoking Proposition~3.5 from \cite{Luan_Yen_2015}, we have 
\begin{equation*}
v \in (({\rm ker}A)^{\perp})^{\perp} \cap \big\{x \in X \mid \langle x^*_i, x \rangle \leq 0, \,  i \in I(x) \big\}.
\end{equation*}
 Since the linear subspace ${\rm ker}A \subset X$ is closed, by \cite[Proposition 2.40]{Bonnans_Shapiro_2000} one gets $(({\rm ker}A)^{\perp})^{\perp}={\rm ker}A$ . Thus, $v \in {\rm ker}A$ and  $\langle x^*_i, v \rangle \leq 0$ for all  $i \in I(x)$. Since $0\in Q_x$, \eqref{separation_Qx} implies that $\langle v^*, v \rangle >0$. Put $x_t=x+tv$, where $t>0$, and note that 
$$Ax_t=Ax+tAv=y.$$
If $i \in I(x)$, then for every $t>0$,
\begin{equation*}
\langle x^*_i, x_t \rangle = \langle x^*_i, x \rangle + t\langle x^*_i, v\rangle = \alpha_i + t\langle x^*_i, v \rangle  \leq \alpha_i. 
\end{equation*}
Since $\langle x^*_j, x \rangle  < \alpha_j$ for any $j \in I\setminus I(x)$, one can find $t>0$ such that  
\begin{equation*}
\langle x^*_j, x_t \rangle = \langle x^*_j, x \rangle + t\langle x^*_j, v \rangle < \alpha_j, \ \, \forall j \in I\setminus I(x).
\end{equation*} 
Hence, for the chosen $t$, we have $x_t \in D$. Since $\langle v^*, x_t-x \rangle=t \langle v^*, v \rangle >0$, it follows that $v^* \notin N_D(x)$. The inclusion $ N_D(x)  \subset  Q_x $ has been proved. Thus $N_D(x)=Q_x$; so $N_D(x)$ is a generalized polyhedral convex cone.  
$\hfill\Box$  

\medskip
During the course of the proof of Theorem~\ref{normal_cone_gpcs}, we have obtained the following result.
\begin{Proposition}\label{rep_normal_cone_gpcs}
	Suppose that $D \subset X$ is a generalized polyhedral convex set given by \eqref{eq_def_gpcs_2}. Then, for any $x \in D$,   
	\begin{equation}\label{eq_rep_normal_cone_gpcs}
	N_D(x)={\rm cone}\{x^*_i \mid i \in I(x)\}+({\rm ker}A)^{\perp}.
	\end{equation}		
\end{Proposition}

\begin{Remark}{\rm In a Banach space setting, formula \eqref{eq_rep_normal_cone_gpcs} has been given in the proof of  \cite[Proposition~3.2]{Ban_Mordukhovich_Song_2011}.	
}\end{Remark}

In connection  with Theorem \ref{normal_cone_gpcs}, one may ask: \textit{If $D \subset X$ is a polyhedral convex set and if  $x \in D$, then $N_D(x)$ is a polyhedral convex cone, or not?} An answer for that question is given in the next statement.

\begin{Proposition}\label{normal_cone_pcs}
	Suppose that $D \subset X$ is a polyhedral convex set and $x \in D$. Then, $N_D(x)$ is a polyhedral convex cone in $X^*$ if and only if $X$ is finite-dimensional.  
\end{Proposition}
\noindent
{\bf Proof.} Since $D$ is a pcs, $D$ can be represented in the form \eqref{eq_def_gpcs_2} with $Y=\{0\}$, $A\equiv 0$, and $y=0$. As ${\rm ker}A=X$ one has $({\rm ker}A)^{\perp}=\{0\}$. So, by Proposition~\ref{rep_normal_cone_gpcs}, $N_D(x)={\rm cone}\{x^*_i \mid i \in I(x)\}.$ Applying Lemma~\ref{rep_pcs_lemma} to $N_D(x)$ in $X^*$, we can assert that $N_D(x)$ is a polyhedral convex cone if and only if the linear subspace $\{0\}$ is of finite codimension in $X^*$, i.e., $X^*$ is finite-dimensional. Since the dual space of any lcHtvs of finite dimension is finite-dimensional \cite[pp.~36--37]{Robertson_Robertson_1964}, we have thus completed the proof. $\hfill\Box$

\medskip
One has the following analogue of Proposition~\ref{rep_normal_cone_gpcs} for polyhedral convex sets.

\begin{Proposition}\label{rep_normal_cone_pcs} 
Suppose that $D \subset X$ is a polyhedral convex set given by \eqref{eq_def_gpcs_2}, where $Y=\{0\}$, $A\equiv 0$, and $y=0$. Then, for every $x \in D$,
	\begin{equation*}\label{eq_rep_normal_cone_pcs}
	N_D(x)={\rm cone}\{x^*_i \mid i \in I(x)\}.	
	\end{equation*}	
\end{Proposition}

\medskip
By the definition of normal cone, we have
\begin{equation}\label{eq_normal_cone_relation}
N_{C_1}(x)+N_{C_2}(x) \subset N_{C_1 \cap C_2}(x),
\end{equation}
for any $x \in C_1 \cap C_2$, where $C_1, C_2$ are convex subsets of $X$. The inclusion \eqref{eq_normal_cone_relation} holds with equality if $X=\mathbb{R}^n$, ${\rm ri}C_1 \cap C_2 \neq \emptyset$ and $C_2$ is polyhedral convex (see \cite[p.~267]{Bertsekas_et_al_2003}), or~$X$ is a lcHtvs and ${\rm int}C_1 \cap C_2 \neq \emptyset$ (see \cite[Proposition~1, p.~205]{Ioffe_Tihomirov_1979}). The next theorem furnishes a property of normal cones to the intersection of gpcs. 

\begin{Theorem}\label{sum_normal_cone_of_gpcs}
	 Let $D_1$ and $D_2$ be two generalized polyhedral convex sets of $X$. For every $x \in D_1 \cap D_2$, 
	\begin{equation}\label{eq_sum_normal_cone_gpcs}
	N_{D_1 \cap D_2}(x)=\overline{N_{D_1}(x)+N_{D_2}(x)}.
	\end{equation}
\end{Theorem}

To prove this result, we need two lemmas.

\begin{Lemma}\label{closed_sum_set}
	Let $C_1, C_2$ be two subsets of a Hausdorff topological vector space $Z$. If~$C_1+ \overline{C_2}$ is closed, then $C_1+ \overline{C_2}=\overline{C_1+C_2}$.
\end{Lemma}
\noindent
{\bf Proof.} Since $C_1+C_2 \subset C_1+ \overline{C_2}$ and since $C_1+ \overline{C_2}$ is a closed set, we see that $\overline{C_1+C_2} \subset C_1+ \overline{C_2}$. Theorem~1.13 (b) from \cite{Rudin_1991} tells us that $\overline{C_1} + \overline{C_2} \subset \overline{C_1+C_2}$; hence $C_1+ \overline{C_2} \subset \overline{C_1+C_2}$. We have thus shown that $C_1+ \overline{C_2}=\overline{C_1+C_2}$. $\hfill\Box$

\begin{Lemma}\label{sum_annihilator}
	Let $M_1, M_2$ be two closed linear subspaces of a locally convex Hausdorff topological vector space $Z$, whose dual space is $Z^*$. Then we have
\begin{equation}\label{eq_sum_annihilator}
(M_1 \cap M_2)^{\perp}=\overline{M_1^{\perp}+M_2^{\perp}}.
\end{equation}	      
\end{Lemma}
\noindent
{\bf Proof.} Applying \cite[formula~(2.32), p.~32]{Bonnans_Shapiro_2000} with  the closed convex cones being replaced by the closed linear subspaces $M_1$ and $M_2$ gives \eqref{eq_sum_annihilator}.  
$\hfill\Box$

\medskip
\noindent
{\bf Proof of Theorem~\ref{sum_normal_cone_of_gpcs}.} For each $k\in \{1,2\}$, since $D_k$ is a gpcs, there exist a continuous linear mapping $A_k$ from $X$ to a lcHtvs $Y_k$, a point $y_k \in Y_k$, a finite index set $I_k$, $x^*_i \in X^*$ and $\alpha_i \in \mathbb{R}$ for $i \in I_k$, such that 
\begin{equation*}
D_k=\{x \in X \mid A_kx=y_k, \langle x^*_i, x \rangle \leq \alpha_i, \, i \in I_k\}.
\end{equation*}     
(We assume that $I_1\cap I_2=\emptyset$.) For $I:=I_1 \cup I_2$,  one has
\begin{equation*}
D_1 \cap D_2=\{x \in X \mid A_kx=y_k, k=1,2, \,  \langle x^*_i, x \rangle \leq \alpha_i, \, i \in I\}.
\end{equation*} 
For each $x \in D_1 \cap D_2$, put $I_k(x)=\{ i \in I_k \mid \langle x^*_i, x \rangle=\alpha_i\}$ for $k=1,2$, $$I(x)=\{ i \in I \mid \langle x^*_i, x \rangle=\alpha_i\},$$ and note that $I(x)=I_1(x) \cup I_2(x)$. On one hand, by  \eqref{eq_rep_normal_cone_gpcs} we have
$$N_{D_1 \cap D_2}(x)={\rm cone}\{x^*_i \mid i \in I(x)\}+({\rm ker}A_1  \cap {\rm ker}A_2)^{\perp}.$$ 
Since $({\rm ker}A_1  \cap {\rm ker}A_2)^{\perp}=\overline{({\rm ker}A_1)^{\perp} + ({\rm ker}A_2)^{\perp}}$ by Lemma~\ref{sum_annihilator}, this implies that
\begin{equation}\label{proof_sum_normal_cone_1}
N_{D_1 \cap D_2}(x)={\rm cone}\{x^*_i \mid i \in I(x)\} + \overline{({\rm ker}A_1)^{\perp} + ({\rm ker}A_2)^{\perp}}.
\end{equation}
On the other hand, applying Proposition \ref{rep_normal_cone_gpcs} to both sets $D_1$ and $D_2$ for $x\in D_1\cap D_2$, we get
\begin{equation}\label{proof_sum_normal_cone_2}
\begin{aligned}
N_{D_1}(x) + N_{D_2}(x)&=\big({\rm cone}\{x^*_i \mid i \in I_1(x)\}+ ({\rm ker}A_1)^{\perp}\big)\\
&\hspace{3cm}+\left({\rm cone}\{x^*_i \mid i \in I_2(x)\}+({\rm ker}A_2)^{\perp}\right)\\
&={\rm cone}\{x^*_i \mid i \in I(x)\} +({\rm ker}A_1)^{\perp} + ({\rm ker}A_2)^{\perp}. 
\end{aligned}
\end{equation}  
Let $C_1:={\rm cone}\{x^*_i \mid i \in I(x)\} $ and $C_2:=({\rm ker}A_1)^{\perp} + ({\rm ker}A_2)^{\perp}$ and observe that $C_1 + \overline{C_2}$ is a generalized polyhedral convex cone in $X^*$ by \cite[Theorem~2.10]{Luan_Yen_2015}. In particular, $C_1 + \overline{C_2}$ is closed. In accordance with Lemma~\ref{closed_sum_set}, $C_1 + \overline{C_2}=\overline{C_1+C_2}$.  
In combination with \eqref{proof_sum_normal_cone_1} and \eqref{proof_sum_normal_cone_2}, this equality justifies \eqref{eq_sum_normal_cone_gpcs}. $\hfill\Box$

\begin{Remark}\label{Rem_sum_normal_cone_gpcs}
		{\rm One may ask: \textit{Whether the closure sign in \eqref{eq_sum_normal_cone_gpcs} can be omitted, or not?} If $X$ is a finite-dimensional space, then $N_{D_1}(x)$ and $N_{D_2}(x)$ are polyhdedral convex cones in the finite-dimensional $X^*$; hence $N_{D_1}(x) + N_{D_2}(x)$ is  polyheral convex by \cite[Corollary~19.3.2]{Rockafellar_1970}. Since $N_{D_1}(x) + N_{D_2}(x)$ is closed,  the closure sign in \eqref{eq_sum_normal_cone_gpcs} is superfluous. However, when $X$ is an infinite-dimensional space, $N_{D_1}(x) + N_{D_2}(x)$ may be non-closed. To see this, one can choose an infinite-dimensional Hilbert space $X$ and two suitable closed linear subspaces $X_1, X_2$ of $X$ so that  $\overline{X_1+X_2}=X$ and $X_1+X_2 \neq X$ (see \cite[Example~3.34]{Bauschke_Combettes_2011} for details). Let $D_i$ be the orthogonal complement of $X_i$, i.e.,  $D_i=\{x\in X\mid \langle x,u\rangle=0\ \forall u\in X_i\}$, for $i=1,2$. It is clear that $D_1, D_2$ are gpcs in $X$ and $D_1 \cap D_2=\{0\}$. Since $N_{D_1}(0)=X_1$ and $N_{D_2}(0)=X_2$, we can assert that $N_{D_1}(0) + N_{D_2}(0)$ is non-closed.}\end{Remark}

In the proof of Theorem~\ref{sum_normal_cone_of_gpcs}, if $D_1$ is a pcs, then we can choose $Y_1=\{0\}$, $A_1\equiv 0$, and $y_1=0$. Since $({\rm ker}A_1)^{\perp}=\{0\}$, one has $\overline{({\rm ker}A_1)^{\perp} + ({\rm ker}A_2)^{\perp}}= ({\rm ker}A_2)^{\perp}$. Hence, \eqref{proof_sum_normal_cone_1} and \eqref{proof_sum_normal_cone_2} imply that 
$N_{D_1 \cap D_2}(x)=N_{D_1}(x) + N_{D_2}(x)$. Thus we have obtained the following result. 
\begin{Theorem}\label{sum_normal_cone_pcs_gpcs}
	Suppose that $D_1 \subset X$ is a polyhedral convex set and $D_2 \subset X$ is a generalized polyhedral convex set. Then, for every $x \in D_1 \cap D_2$, 
		\begin{equation*}\label{eq_sum_normal_cone_pcs_gpcs}
		N_{D_1 \cap D_2}(x)=N_{D_1}(x)+N_{D_2}(x).
		\end{equation*}
\end{Theorem}	

\medskip
Following  \cite[p.~34]{Robertson_Robertson_1964}, we define the \textit{polar} of a nonempty set $C$ by 
$$C^o:=\big\{x^* \in X^* \mid \langle x^*, x \rangle \leq 1, \, \forall x \in C \big\}.$$
Evidently, $C^o$ is a weakly$^*$-closed convex set containing $0$. If $C$ is a cone, then one has $C^o=N_C(0)$. 

\medskip The forthcoming proposition extends \cite[Corollary~19.2.2]{Rockafellar_1970} to a lcHtvs setting.

\begin{Proposition}\label{polar_set_of_gpcs}
	The polar of a nonempty generalized polyhedral convex set is a generalized polyhedral convex set. 
\end{Proposition}
\noindent
{\bf Proof.} Suppose that $D \subset X$ is given by \eqref{rep_D}. Then we have
\begin{equation}\label{eq_polar_gpcs}
D^o=\big\{x^* \in X_0^{\perp} \mid \langle x^*, u_i \rangle \leq 1, \, i=1,\dots,k, \, \langle x^*, v_j \rangle \leq 0, \, j=1,\dots,\ell \big\}.
\end{equation}
Indeed, take any $x^* \in D^o$. The inequalities $\langle x^*, u_i \rangle \leq 1$, $i=1,\dots,k$, are valid, because $u_i  \in D$. As the linear functional $\langle x^*, \cdot\rangle $ is bounded from the above on $D$, Theorem~3.3 from \cite{Luan_Yen_2015} shows that the linear programming problem
$$\max\{\langle x^*,x \rangle \mid  x \in D\}$$ has a solution. Therefore,
by \cite[ Proposition~3.5]{Luan_Yen_2015}, we have $x^* \in X_0^{\perp}$ and $\langle x^*, v_j \rangle \leq 0$ for all $j=1,\dots,\ell$. The inclusion ``$\subset$'' in \eqref{eq_polar_gpcs} has been proved. To obtain the opposite inclusion, take any $x^*$ from the set on the right-hand side of \eqref{eq_polar_gpcs}. By~\eqref{rep_D}, for each $x \in D$, there exist  nonnegative numbers $\lambda_1,\dots, \lambda_k,$ $\mu_1,\dots, \mu_\ell$, and a vector $x_0 \in X_0$, such that $\sum\limits_{i=1}^k \lambda_i=1$ and
$x=\sum\limits_{i=1}^k \lambda_i u_i + \sum\limits_{j=1}^{\ell} \mu_j v_j +x_0$. Since
\begin{equation*}
\begin{aligned}
\langle x^*, x \rangle&= \sum\limits_{i=1}^k \lambda_i \langle x^*, u_i \rangle + \sum\limits_{j=1}^{\ell} \mu_j \langle x^*, v_j \rangle + \langle x^*, x_0 \rangle \\
&\leq \sum\limits_{i=1}^k \lambda_i \langle x^*, u_i \rangle \leq \sum\limits_{i=1}^k \lambda_i =1,
\end{aligned}
\end{equation*}
we see that $\langle x^*, x \rangle \leq 1$ for any $x \in D$; hence $x^* \in D^o$. This completes the proof of~\eqref{eq_polar_gpcs}. The fact that $D^o$ is a gpcs in $X^*$ follows from \eqref{eq_polar_gpcs} and Definition \ref{Def_gpcs}.  
$\hfill\Box$

\medskip
According to \cite[p.~172]{Ioffe_Tihomirov_1979}, the {\it conjugate function} (or the \textit{Young-Fenchel transform function}) of a function $f:X\to \overline{\mathbb{R}}$ is the function $f^{*}: X^* \to \overline{\mathbb{R}}$ given by 
\begin{equation*}\label{eq_def_conjugate_fun_1}
f^{*}(x^*)=\sup\big\{ \langle x^*, x \rangle - f(x) \mid x \in X \big\}.
\end{equation*} 
It is well known  \cite[Proposition~3, p.~174]{Ioffe_Tihomirov_1979} that if $f$ is proper convex and lower semicontinuous (i.e., ${\rm epi}f$ is a closed set), then $f^{*}$ is also a proper convex lower semicontinuous function. 
It is clear that $
f^{*}(x^*)=\sup\big\{ \langle x^*, x \rangle - f(x) \mid x \in {\rm dom}f \big\}$ for any $x^* \in X^*$. 

\begin{Theorem}\label{conjugate_fun} The conjugate function of a proper generalized polyhedral convex function is a proper generalized polyhedral convex function.   
\end{Theorem}
\noindent
{\bf Proof.} Suppose that  $f:X\to \overline{\mathbb{R}}$ is a proper gpcf. Then $f^*$ is a proper convex function. Moreover, due to Theorem~\ref{rep_gpcf_gpl}, $f$ is a generalized piecewise linear function. So, there exist nonempty gpcs $D_1, \dots, D_m$ in $X$, $v^*_k \in X^*$, $\beta_k \in \mathbb{R}$, $k=1,\dots,m$, such that  ${\rm dom}f=\bigcup\limits_{k=1}^m D_k$ and $f(x)=\langle v^*_k, x \rangle + \beta_k$ for every $x \in D_k$, $k=1,\dots, m$. For each~$k$, by \cite[Theorem 2.7]{Luan_Yen_2015} we can find finite index sets $I_k$ and $J_k$, points $u_i\in X$ with $i \in I_k$, vectors $v_j \in X$ with $j \in J_k$, and a closed linear subspace $X_{0,k}$ in $X$, such that
$$D_k={\rm conv}\{u_i  \mid i \in I_k \} + {\rm cone}\{v_j  \mid j \in J_k \} + X_{0,k}.$$
(We assume that $I_k\cap I_{\ell}=\emptyset$ and $J_k\cap J_{\ell}=\emptyset$ whenever $k \neq \ell$.)
For every $k\in \{1,\dots,m\}$, consider the function $\varphi_k(x^*)=\sup\big\{ \langle x^*, x \rangle - \langle v^*_k, x \rangle - \beta_k \mid x \in  D_k\big\}$ defined on $X^*$ and observe that $x^* \in {\rm  dom}\varphi_k$ if and only if the linear functional $\langle x^*, \cdot \rangle - \langle v^*_k, \cdot \rangle - \beta_k$ is bounded from the above on $D_k$. The latter is equivalent to the property that the linear programming problem
$\max\{\langle x^* - v^*_k, x \rangle - \beta_k \mid  x \in D_k\}$ has a solution (see \cite[Theorem~3.3]{Luan_Yen_2015}). Therefore,
by \cite[Proposition~3.5]{Luan_Yen_2015} we get
\begin{equation*}
{\rm dom}\varphi_k=\left\{x^* \in X^* \mid x^*-v^*_k \in X_{0, k}^{\perp}, \, \langle x^*-v^*_k, v_j \rangle \leq 0,  \, j \in J_k   \right\}.
\end{equation*} 
As $f^{*}(\cdot)=\max\big\{ \varphi_k(\cdot) \mid k=1,\dots,m \big\}$, one has ${\rm dom}f^*=\bigcap\limits_{k=1}^m {\rm dom} \varphi_k$; hence
\begin{equation*}
{\rm dom}f^*=\left\{x^* \in \bigcap_{k=1}^m \big(v^*_k+X_{0, k}^{\perp} \big)\mid  \langle x^*,v_j\rangle \leq \langle  v_k^*,v_j\rangle, k=1,\dots,m, \, j \in J_k \right\}.
\end{equation*} 
Since $\bigcap\limits_{k=1}^m \big(v^*_k+X_{0, k}^{\perp} \big)$ is a closed affine subspace of $X^*$, we can assert that ${\rm dom}f^*$ is a gpcs. For every $x^* \in {\rm dom}\varphi_k$, it is a plain matter to show that $$\varphi_k(x^*)=\max\{ \langle x^*-v^*_k, u_i \rangle -\beta_k \mid i \in I_k\}.$$ Therefore,
		\begin{equation}\label{rep_conjugate_func}
		f^*(x^*)=\begin{cases}
		\max \left\{ \langle x^*, u_i \rangle -f(u_i) \mid k=1,\dots,m, \, i \in I_k\right\}  &\text{if } x^* \in {\rm dom}f^*,\\
		+\infty & \text{if } x^* \notin {\rm dom}f^*.
		\end{cases}
		\end{equation}	  
Since $f^*$ is a proper function with ${\rm dom}f^*$ being a gpcs, using Theorem~\ref{rep_gpcf} and \eqref{rep_conjugate_func} we can conclude that $f^*$ is a gpcf. 
$\hfill\Box$

\begin{Remark}{\rm
	Theorem~\ref{conjugate_fun}  is a generalization of Theorem~19.2 from \cite{Rockafellar_1970}, where the case $X=\mathbb{R}^n$ was treated.
}\end{Remark}

In the remaining part of this section, we will study subdifferentials of  generalized polyhedral convex functions. It is well known that the subdifferential of a convex function is the basis for optimality conditions and  other issues in convex programming. On account of \cite[p.~46]{Ioffe_Tihomirov_1979}, a linear functional $x^* \in X^*$ is said to be a \textit{subgradient}  of a proper convex function $f$ at $x \in {\rm dom}f$ if 
$$ \langle x^*, u - x \rangle \leq f(u)-f(x) \quad (u \in X).$$ This condition is equivalent to the simple geometric property that the graph of the affine function $h(u)= f(x)+\langle x^*, u - x \rangle$ forms a non-vertical supporting hyperplane to ${\rm epi}f$ at the point $(x, f(x))$; see \cite[pp.~214--215]{Rockafellar_1970}. The \textit{subdifferential} of $f$ at $x$, denoted by $\partial f(x)$, is the set of all the subgradients of $f$ at $x$. From the defintion it follows that $\partial f(x)$ is a weakly$^*$-closed convex set (see \cite[p.~81]{Bonnans_Shapiro_2000}). Moreover, by \cite[Propostion~1, p.~197]{Ioffe_Tihomirov_1979}, $x^* \in \partial f(x)$ if and only if $f(x)+f^*(x^*)=\langle x^*, x \rangle$. If $C$ is a nonempty convex subset of~$X$ then, for any $x \in C$, one has $\partial \delta(x, C)=N_C(x)$, where~$\delta(\cdot, C)$ is the indicator function of $C$.

\medskip
Based on Theorem \ref{rep_gpcf}, the next theorem provides us with a formula for the subdifferential of a gpcf.
\begin{Theorem}\label{rep_subdifferential_gpcf} Suppose that $f$ is a proper generalized polyhedral convex function with ${\rm dom}f= \left\{ x \in X \mid Ax=y, \, \langle x^*_i, x \rangle \leq \alpha_i, \, i=1,\dots,p\right\}$ and
	\begin{equation*}
	f(x)=	\max\left\{ \langle v_j^*, x \rangle +\beta_j\mid j=1,\dots,m \right\}\quad (x \in {\rm dom}f),
	\end{equation*} 
where $A$ is a  continuous linear mapping from $X$ to a locally convex Hausdorff topological vector space $Y$, $y \in Y$, $x^*_i \in X^*, \alpha_i \in \mathbb{R}$, $i=1,\dots,p$, $v_j^* \in X^*, \beta_j \in \mathbb{R}$, $j=1,\dots,m$. Then, for every $x \in {\rm dom}f$,  
\begin{equation}\label{eq_rep_subdifferential_gpcf}
\partial f(x)={\rm conv}\{v_j^* \mid  j \in J(x)\}+{\rm cone}\{x^*_i \mid  i \in I(x)\}+({\rm ker}A)^{\perp},
\end{equation}	
where $I(x)=\big\{i \in \{1,\dots,p\} \mid \langle x_i^*, x \rangle= \alpha_i\big\}$ and $$J(x)=\big\{j \in \{1,\dots,m\} \mid \langle v_j^*, x \rangle +\beta_j = f(x)\big\}.$$ In particular, if $Y=\{0\}$, $A \equiv 0$ and $y=0$ (the case where ${\rm dom}f$ is a polyhedral convex set) then, for any $x \in {\rm dom}f$,  
\begin{equation*}\label{eq_rep_subdifferential_pcf}
\partial f(x)={\rm conv}\{v_j^* \mid  j \in J(x)\}+{\rm cone}\{x^*_i \mid  i \in I(x)\}.
\end{equation*}	  
\end{Theorem}  
\noindent
{\bf Proof.} Fix any $x \in D:={\rm dom}f$. As the functions $f_j(\cdot):=\langle v_j^*, \cdot \rangle +\beta_j$ is continuous at~$x$ for all $j=1,\dots,m$, the function $\widetilde{f}(\cdot):=\max\left\{ f_j(\cdot)\mid j=1,\dots,m \right\}$ is also continuous at $x$. Hence,  applying the Moreau--Rockafellar theorem \cite[Theorem~1, p.~200]{Ioffe_Tihomirov_1979} to the sum  $\widetilde{f}(\cdot)+\delta(\cdot,D)=f(\cdot)$,  one gets
\begin{equation}\label{proof_rep_subdifferential_1}
\partial f(x) = \partial \widetilde{f}(x) +  \partial \delta(x,D)= \partial \widetilde{f}(x)  + N_D(x).
\end{equation}
Since $\partial f_j(\cdot) \equiv \{v_j^*\}$ for all $j=1,\dots,m$, by \cite[Theorem 3, p.~201]{Ioffe_Tihomirov_1979} we obtain
\begin{equation*}
\partial \widetilde{f}(x) = \overline{\rm conv} \Big(\bigcup\limits_{j \in J(x)} \partial f_j(x) \Big)=\overline{\rm conv} \big\{ v_j^* \mid j \in J(x) \big\}.
\end{equation*}
 On one hand, according to \cite[Theorem~2.7]{Luan_Yen_2015}, ${\rm conv} \big\{ v_j^* \mid j \in J(x) \big\}$ is a gpcs (hence it is closed). So, $ \partial \widetilde{f}(x) = {\rm conv} \big\{ v_j^* \mid j \in J(x) \big\}.$ On the other hand, in accordance with Proposition~\ref{rep_normal_cone_gpcs},
 $ N_D(x)={\rm cone}\{x^*_i \mid i \in I(x)\}+({\rm ker}A)^{\perp}.$ Therefore, from \eqref{proof_rep_subdifferential_1} one obtains \eqref{eq_rep_subdifferential_gpcf}.  $\hfill\Box$

\medskip
From \eqref{eq_rep_subdifferential_gpcf} and \cite[Theorem 2.7]{Luan_Yen_2015} it follows that $\partial f(x)$ is a gpcs in $X^*$. Thus we have proved the following result, which is a known one \cite[Theorem~23.10]{Rockafellar_1970} in the case where $X=\mathbb{R}^n$.

\begin{Proposition}
	If $f$ is a proper generalized polyhedral convex function on $X$ and if $x \in {\rm dom}f$, then $\partial f(x)$ is a generalized polyhedral convex set.
\end{Proposition}

By the definition of subdifferential, if $f_1, \dots, f_m$ are proper convex functions on~$X$ then, for every $x \in \bigcap\limits_{i=1}^m {\rm dom}f_i$,
\begin{equation}\label{eq_subdifferential_relation}
\partial f_1(x)+ \dots+ \partial f_m(x) \subset	\partial (f_1+\dots+f_m)(x).
\end{equation} Since the set on the right-hand side of \eqref{eq_subdifferential_relation} is weakly$^*$-closed, one has 
\begin{equation*}
\overline{\partial f_1(x)+ \dots+ \partial f_m(x)} \subset	\partial (f_1+\dots+f_m)(x).
\end{equation*}
The above-cited Moreau--Rockafellar theorem tells us that~\eqref{eq_subdifferential_relation} holds with equality if there exists $x_0 \in \bigcap\limits_{i=1}^m {\rm dom}f_i$ such that all the functions $f_1, \dots, f_m$ except, possibly, one are continuous at $x_0$. The specific structure of  generalized polyhedral convex functions allows one to have a subdifferential sum rule without the continuity assumption.

\begin{Theorem}\label{sum_subdifferentials_gpcf}
Let $f_1, \dots, f_m$ be proper generalized polyhedral convex functions. Then, for any $x \in \bigcap\limits_{i=1}^m {\rm dom}f_i$, 
	\begin{equation}\label{eq_sum_subdifferentials_gpcf}
	\partial (f_1+f_2+\dots+f_m)(x)=\overline{\partial f_1(x)+ \partial f_2(x)+\dots+ \partial f_m(x)}.
	\end{equation}
\end{Theorem}
\noindent
{\bf Proof.} According to Theorem \ref{sum_two_gcpf}, the sum of two proper generalized polyhedral convex functions whose effective domains have at least one common point, is again a proper generalized polyhedral convex function. So, to obtain the desired result, it suffices to prove \eqref{eq_sum_subdifferentials_gpcf} for $m=2$ and then proceed by induction.  

For each $i=1,2$, by Theorem~\ref{rep_gpcf}, $D_i:={\rm dom}f_i$ is a gpcs and there exist $v^*_{i,j} \in X^*$, $\beta_{i,j} \in \mathbb{R}$, $j=1,\dots,k_i$, such that $f_i(x)=\widetilde{f_i}(x)+\delta(x, D_i)$, where
$$\widetilde{f_i}(x)=\max\left\{ \langle v^*_{i,j}, x \rangle +\beta_{i,j} \mid j=1,\dots,k_i \right\} \quad (x \in X).$$
Clearly, $\widetilde{f_1}$ and $\widetilde{f_2}$ are proper and continuous on $X$. Let $x \in D_1 \cap D_2$. On one hand, by the formula $f_1+f_2=\widetilde{f_1}+\widetilde{f_2}+\delta(\cdot, D_1 \cap D_2)$ and by the Moreau--Rockafellar theorem, 
\begin{equation}\label{eq_sum_subdifferentials_gpcf_1}
\begin{aligned}
\partial (f_1+f_2)(x) &= \partial \widetilde{f_1}(x) +  \partial \widetilde{f_2}(x)+ \partial \delta(x,D_1 \cap D_2)\\
&= \partial \widetilde{f_1}(x) +  \partial \widetilde{f_2}(x)+ N_{D_1 \cap D_2}(x).
\end{aligned}
\end{equation}
Since $N_{D_1 \cap D_2}(x)=\overline{N_{D_1}(x)+N_{D_2}(x)}$ by Theorem~\ref{sum_normal_cone_of_gpcs}, this implies that
\begin{equation}\label{eq_sum_subdifferentials_gpcf_2}
\partial (f_1+f_2)(x) = \partial \widetilde{f_1}(x) +  \partial \widetilde{f_2}(x)+ \overline{N_{D_1}(x)+N_{D_2}(x)}.
\end{equation}
On the other hand, applying the Moreau--Rockafellar theorem to the proper convex functions $f_1=\widetilde{f_1}+\delta(\cdot, D_1)$  and $f_2=\widetilde{f_2}+\delta(\cdot, D_2)$, we obtain
\begin{equation}\label{eq_sum_subdifferentials_gpcf_3}
\partial f_1(x)+\partial f_2(x)=\partial \widetilde{f_1}(x) + N_{D_1}(x)+ \partial \widetilde{f_2}(x)+N_{D_2}(x).
\end{equation}
Put $C_1=\partial \widetilde{f_1}(x) + \partial \widetilde{f_2}(x)$ and $C_2=N_{D_1}(x)+N_{D_2}(x)$. From \eqref{eq_sum_subdifferentials_gpcf_2} and the closedness of $\partial (f_1+f_2)(x)$, it follows that $C_1 + \overline{C_2}$ is closed. Then, according to Lemma~\ref{closed_sum_set}, $C_1 + \overline{C_2}=\overline{C_1+C_2}$. Combining this equality with \eqref{eq_sum_subdifferentials_gpcf_2} and \eqref{eq_sum_subdifferentials_gpcf_3} yields \eqref{eq_sum_subdifferentials_gpcf}. $\hfill\Box$

\medskip
In the last proof, if $f_1$ is a pcf, then $D_1$ is a pcs by virtue of Theorem \ref{rep_gpcf}. Hence, by Theorem~\ref{sum_normal_cone_pcs_gpcs} we have $N_{D_1 \cap D_2}(x)=N_{D_1}(x)+N_{D_2}(x)$. Therefore, using  \eqref{eq_sum_subdifferentials_gpcf_1} and~\eqref{eq_sum_subdifferentials_gpcf_3}, we can obtain formula \eqref{eq_sum_subdifferentials_gpcf} in the case $m=2$ with no closure sign on the right-hand side. Thus, the following result is valid.

\begin{Theorem}\label{sum_subdifferentials_pcf_gpcf}
	Suppose that $f_1$ is a proper polyhedral convex function and $f_2$ is a proper generalized polyhedral convex function. Then, for any $x \in ({\rm dom}f_1) \cap ({\rm dom}f_2)$, 
	\begin{equation*}\label{eq_sum_subdifferentials_pcf_gpcf}
	\partial (f_1+f_2)(x)=\partial f_1(x)+ \partial f_2(x).
	\end{equation*}
\end{Theorem}

\section*{Acknowledgements}
This work was supported by National Foundation for Science $\&$ Technology Development (Vietnam) and China Medical University (Taichung, Taiwan). The second author was partially supported by the Grant MOST 105-2115-M-039-002-MY3.

\end{document}